\documentclass[12pt]{article}
\usepackage{amsthm}
\usepackage{amsmath}
\usepackage{amsfonts}
\usepackage{amssymb}
\usepackage{url}
\usepackage{epsfig}

\newtheorem{thm}{Theorem}[section]
\newtheorem{cor}[thm]{Corollary}
\newtheorem{lem}[thm]{Lemma}
\newtheorem{prop}[thm]{Proposition}

\theoremstyle{definition}
\newtheorem{example}[thm]{Example}
\newtheorem{remark}[thm]{Remark}

\newtheorem{problem}[thm]{Problem}

\newcommand{\R}{\mathbb{R}}

\newcommand{\Z}{\mathbb{Z}}
\newcommand{\N}{\mathbb{N}}
\newcommand{\Q}{\mathbb{Q}}
\newcommand{\PP}{\mathbb{P}}

\date{}

\begin{document}
\bibliographystyle{plain}
\title{\bf Sagbi Bases of Cox-Nagata Rings}
\author{Bernd Sturmfels \ and \ Zhiqiang Xu}

\maketitle

\begin{abstract}
We degenerate Cox-Nagata rings to toric algebras by means of
sagbi bases  induced by configurations over the rational function field.
For del Pezzo surfaces, this degeneration implies
the Batyrev-Popov conjecture
that these rings are presented by ideals of quadrics.
For the blow-up of projective $n$-space at $n+3$ points,
sagbi bases of Cox-Nagata rings establish a link between
the Verlinde formula and phylogenetic
algebraic geometry, and we use this to answer questions
due to D'Cruz-Iarobbino and
Buczy\'nska-Wi\'sniewski.
Inspired by the zonotopal algebras of Holtz and Ron,
our study emphasizes explicit computations, and
offers a new approach to Hilbert functions of fat points.
\end{abstract}

\section{Powers of linear forms}

We fix $n$ vector fields on a $d$-dimensional space
with coordinates $(z_1,\ldots,z_d)$:
$$ \ell_j \,\,= \,\, \sum_{i=1}^d a_{ij} \frac{\partial}{\partial z_i}
\qquad \hbox{for} \,\, j = 1,2,\ldots,n. $$
Here the coefficients $a_{ij}$ are scalars in a field $K$,
which we assume to have characteristic zero.
We regard $\ell_1,\ell_2, \ldots, \ell_n$ as
linear forms in the polynomial ring $S = K[\partial_1, \ldots, \partial_d]$,
where $\,\partial_i \,=\,\partial/\partial z_i$.
We further assume that  $\{\ell_1,\ldots,\ell_n\}$ spans the space $S_1$ of all linear forms in $S$.
This implies $d \leq n$.
We are interested in the family of zero-dimensional polynomial ideals
$$ I_u \,\, := \,\,
\bigl\langle  \ell_1^{u_1+1},   \ell_2^{u_2+1}, \ldots,
\ell_n^{u_n+1} \bigr\rangle \quad \subset \,\, S,$$ where $u =
(u_1,u_2,\ldots,u_n)$ runs over the set $\N^n$ of non-negative
integer vectors. The ideal  $I_u$ represents a system of linear partial differential
equations with constant coefficients, and we consider its solution space
$$ (I_u)^\perp \quad = \quad \bigoplus_{r=0}^\infty I_{(r,u)}^\perp
\,\,\subset \,\, K[z_1,\ldots,z_d] \,\, =: \,\, K[z]. $$
Here $I_{(r,u)}^\perp$ is the space of all polynomials
 that are homogeneous
of degree $r$ and are annihilated by $u_i+1$ repeated applications of the
vector field $\ell_i$. Our principal object of interest is the associated
{\em counting function}
\begin{equation}
\label{countingfunction} \psi \,: \, \N^{n+1} \rightarrow \N \,,\,\,
(r,u) \mapsto {\rm dim}_K I_{(r,u)}^\perp,
\end{equation}
where $\N$ denotes the set of non-negative integers. Note that
$\psi(r,u)$ is the value of the Hilbert function of $S/I_u$ at the
integer~$r$. Here are some examples.

\begin{example}
\label{ex:eins} ($n=d$) If the given linear forms are linearly
independent then we can assume $\ell_1 = z_1,\ldots,\ell_d = z_d$,
so that $I_{(r,u)}^\perp$ is spanned by all monomials $z^v = z_1^{v_1}
\cdots z_d^{v_d}$ whose degree $\,v_1 + \cdots + v_d \,$ equals
$r$. The number $\psi(r,u)$ of these monomials equals the
coefficient of $q^r$ in  the generating function
$$ \frac{(1-q^{u_1+1})(1-q^{u_2+1}) \cdots (1-q^{u_d+1})}{(1-q)^d}. $$
Note that $\psi : \N^{d+1} \rightarrow \N$ is a piecewise polynomial
function of degree $d - 1$. \qed
\end{example}

\begin{example}
\label{ex:zwei} ($d=2$) For pairwise linearly independent binary
linear forms,
\begin{equation}
\label{plusplus}
 \psi(r,u) \quad = \quad \bigl( \,r+1 - \sum_{i=1}^n(r-u_i)_+ \, \bigr)_+  \, ,
 \end{equation}
where $(x)_+ = {\rm max}\{x,0\}$, so $\psi : \N^{n+1} \rightarrow
\N$ is a piecewise linear function. \qed
\end{example}

The formulas in Examples \ref{ex:eins} and \ref{ex:zwei} are
well-known. They can be found, for instance, in the article on our problem
by Kuttler and Wallach \cite{KW}. The next example illustrates the
new type of formulas  to be derived in this paper.

\begin{example}
\label{ex:drei} ($d=3, n=5$) For five general linear forms in
$K[\partial_1,\partial_2,\partial_3]$, the piecewise quadratic function $\psi : \N^{6}
\rightarrow \N$  is given by counting lattice points in the
following convex polygon.  The value $\psi(r,u_1,u_2,\ldots,u_5)$
equals
\begin{eqnarray*}
 \# \bigl\{ \, (x,y) \in  \Z^2 \, \, : & \!\!\!\!
 {\rm max}(0, 2r-u_1-u_2-u_3) \,\leq \,x\,\leq {\rm min}(u_4,u_5) \quad \hbox{and} \\ &
 \!\!\!\!
  {\rm max}(0, 2r-u_3-u_4-u_5) \,\leq \,y\, \leq\, {\rm min}(u_1,u_2) \quad \hbox{and}\\ &
 \!\!\!\!\!\!\!\!\!\!\!\!\!\!\!\!\!
r-u_3 \,\leq\, x+y \,\leq \, r \,\,\, \hbox{and} \,\,\,
 r-u_1-u_2 \,\leq\, y-x\, \leq  \,u_4 + u_5 -r \bigr\} .
\end{eqnarray*}
The derivation of this Ehrhart-type formula is presented in Section~4.
 \qed
\end{example}

We now discuss the organization of this paper and we summarize its
main contributions. In Section 2 we introduce the Cox-Nagata ring,
and we express $\psi$ as the multigraded Hilbert function of that
ring. For $d \geq 3$, this is the Cox ring of the variety gotten from
$\PP^{d-1}$ by blowing up the $n$ points dual to the linear forms
$\ell_i$, and we review some relevant material on fat points
with a fixed support, and on Weyl groups generated by Cremona
transformations.

In Section 3 we present our technique, which is to examine the flat
family induced by taking linear forms defined over $K = \Q(t)$.
Under favorable circumstances, this degeneration has desirable
combinatorial properties, and the minimal generators of the
Cox-Nagata ring form a sagbi basis. In Sections 4, 5 and 6 we demonstrate
that this happens for del Pezzo surfaces, which arise by blowing
up $n \leq 8$ general points in $\PP^2$. We present a proof,
found independently from that given in \cite{TVV}, of the
Batyrev-Popov conjecture \cite{BP, LV, Vera} which states
that the presentation ideals of these Cox-Nagata rings are quadratically generated.
Explicit formulas like that in Example \ref{ex:zwei} are obtained for
the number of sections of line bundles on all del Pezzo surfaces.

In Section 7 we treat the case of  $n=d+2$ points in $\PP^{d-1}$.
We resolve a problem left open by  Buczy\'nska and Wi\'sniewski in \cite{BW},
by showing that the phylogenetic varieties in \cite{BW, SS} arise as flat limits of the
Cox-Nagata rings constructed by Castravet and Tevelev in \cite{CT}.
We also prove a conjecture of D'Cruz and Iarobbino \cite{DI} on Hilbert
functions of fat points. A key tool is an interpretation of the
Verlinde formula \cite{Mukbook} suggested to us by Jenia Tevelev.

In Section 8 we turn to the original motivation which started this project,
namely, the work on zonotopal algebra by Holtz and Ron \cite{HR}.
In their setting, the linear forms $\ell_i$ represent all the hyperplanes that are spanned by
a configuration of vectors, and when $\psi$ is restricted to a certain
linear subspace, matroid theory yields a
particularly beautiful formula for its values.

\section{Cox-Nagata rings}

Let $G$ be the space of linear relations on the linear forms
$\ell_i$. Thus $G$ is the subspace of $K^n$ which consists of all
vectors $\lambda  = (\lambda_1, \lambda_2, \ldots, \lambda_n)$ that
satisfy
$$  \lambda_1 \ell_1(z) + \lambda_2 \ell_2(z) + \cdots +  \lambda_n \ell_n(z)  \,\,\, = \,\,\,   0 . $$
The $K$-vector space $G$ has dimension $ n-d$. We regard it is
an additive group. This  group acts on the polynomial ring in $2n$
variables,
$\, R \, = \,\, K [x_1,x_2,\ldots,x_n, y_1,y_2,\ldots,y_n]$,
by the following {\em Nagata action} (cf.~\cite{CT,  Mukai, Nag}):
$$ x_i \mapsto x_i \quad \hbox{and} \quad
y_i \mapsto y_i + \lambda_i x_i \quad \hbox{for all} \,\,\lambda \in
G .$$ Let $R^G$ denote the subring of $R$ consisting of all
polynomials that are fixed by this action. The {\em Cox-Nagata ring}
$R^G$ is a multigraded ring with respect to the grading which is
induced by  following $\Z^{n+1}$-grading on $R$:
$${\rm deg}(x_i) \, = \, e_i \,\,\,{\rm and} \,\,\, {\rm deg}(y_i) \, = \, e_0 + e_i
\quad \hbox{for} \,\, i = 1,2,\ldots,n. $$ Here $e_0,e_1,\ldots,e_n$
denotes the standard basis of $\Z^{n+1}$, and $R^G_{(r,u)}$ is the
finite-dimensional space of invariant polynomials of multidegree $r
e_0 + \sum_{i=1}^n u_i e_i$. We start out by presenting an
elementary proof of the following  result.

\begin{prop}
There exists a natural isomorphism of $K$-vector spaces
\begin{equation}
\label{iso}
 I_{(r,u)}^\perp \quad \simeq  \quad R^G_{(r,u)}.
 \end{equation}
Hence $ \psi$ is the $\Z^{n+1}$-graded Hilbert function of the
Cox-Nagata ring $R^G$.
\end{prop}

\begin{proof}
We introduce the polynomial ring $K[Y] = K[Y_1, Y_2,\ldots,Y_n]$ and
we let $L_G$ denote the ideal generated by all linear forms
$\lambda_1 Y_1 + \cdots +  \lambda_n Y_n $, where $\lambda$ runs
over $G$. Then $K[Y] /L_G$ is isomorphic to $S = K[z]$ under the
$K$-algebra homomorphism which sends $Y_i$ to $\ell_i(z)$, and this
induces an isomorphism
$$ S/I_u  \,\,\,\simeq \,\,\, K[Y]/\bigl(L_G + \langle Y_1^{u_1+1}, Y_2^{u_2+1}, \ldots, Y_n^{u_n+1}
\rangle\bigr). $$ Therefore the solution space $I_{(r,u)}^\perp$ is
isomorphic to the space of homogeneous polynomial $f$ of degree $r$
in $K[Y]$ that are invariant under the $G$-action $\,Y_i \mapsto Y_i
+ \lambda_i\,$ and that satisfy ${\rm deg}_{Y_i}(f) \leq u_i$ for $i
= 1,2,\ldots,n$. The map
$$ f \,\,\mapsto \,\, f \bigl( \frac{y_1}{x_1} ,\ldots,\frac{y_n}{x_n} \bigr) \cdot
x_1^{u_1} x_2^{u_2} \cdots x_n^{u_n} $$ defines an isomorphism
to the component $R^G_{(r,u)}$ of the Cox-Nagata ring.
\end{proof}

The isomorphism (\ref{iso}) has the following explicit description.
Let $ (a_{ij})$ be the $d \times n$-matrix over $K$ such that
$\,\ell_j = \sum_{i=1}^d a_{ij} \partial_i\,$ for $j=1,\ldots,d$. Then
\begin{equation}
\label{iso2} g(\partial_1,\ldots, \partial_d) \,\,\mapsto \,\, g \bigl(
\sum_{j=1}^n a_{1j} \frac{y_j}{x_j},\ldots,
 \sum_{j=1}^n a_{dj} \frac{y_j}{x_j}\bigr)
 x_1^{u_1} x_2^{u_2} \cdots x_n^{u_n} .
\end{equation}
takes the solution space $\, I_{(r,u)}^\perp \,$ bijectively onto
the graded component $\, R^G_{(r,u)}$. This map can be inverted
precisely for those elements of $R$ that lie in the invariant ring
$R^G$. In this manner, every choice of a homogeneous $K$-basis of
the $\N^{n+1}$-graded $K$-algebra $R^G$ determines a basis of $\,
I_{(r,u)}^\perp \,$ for all $r,u$.

We now explain why we chose the name {\em Cox-Nagata ring} for
$R^G$. In 1959 Nagata resolved Hilbert's $14^{\rm th}$ problem
whether the ring of polynomial invariants of any matrix group is
finitely generated \cite{Nag}. He showed that $R^G$ is not finitely
generated when $G$ is a generic subspace of $K^d$ and $ d =
3, n=16$. The final and precise statement along these lines is due
to Mukai \cite{Mukai}.

\begin{thm} {\rm \cite{CT, Mukai} } Suppose that $G$ is a generic subspace of codimension
$d$ in $K^n$. Then the Cox-Nagata ring $R^G$ is finitely generated
if and only if
\begin{equation}
\label{fracineq}
 \frac{1}{2} \,\, + \,\,\frac{1}{d} \,\,+ \,\,\frac{1}{n-d} \quad > \quad 1 .
 \end{equation}
\end{thm}

The name of David Cox is attached to the ring $R^G$ because of the
following geometric interpretation. Let $\PP^{d-1}$ denote the
projective space whose points are equivalent classes of linear forms
in $S = K[\partial_1,\ldots,\partial_d]$ modulo scaling. Our linear form $\,\ell_j
= \sum_{i=1}^n a_{ij} \partial_i\,$ is represented by the point $(a_{1j}
:a_{2j}: \cdots : a_{dj})$ in $ \PP^{d-1}$. Let $P_j $ denote the
homogeneous prime ideal in $S$ of that point, that is, $P_j$ is the
ideal generated by the $2 \times 2$-minors of the $2 \times
n$-matrix
$$ \begin{pmatrix} \partial_1 & \partial_2 & \cdots & \partial_d \\
 a_{1j} & a_{2j} & \cdots & a_{dj}
 \end{pmatrix} .$$
The following classical fact relates our problem to the study of
ideals of fat points \cite{DI}. We learned Lemma~\ref{ensalem} from
Ensalem and Iarrobino~\cite[Thm.~1]{EI}.

\begin{lem} \label{ensalem} The polynomial solutions of degree $r$ to the linear
partial differential equations with constant coefficients expressed
by $I_u$ are precisely those polynomials that vanish of order $\geq
r-u_j$ at the point $\ell_j$ for all $j$. In symbols,
 \begin{equation}
 \label{sectionspace1}
  I_{(r,u)}^\perp \quad = \quad
 \bigl( P_1^{r-u_1} \,\cap \, P_2^{r-u_2} \, \cap \,\cdots \,\cap \,P_n^{r-u_n}\bigr)_r .
\end{equation}
\end{lem}

Suppose $d \geq 3$ and let $X_G$ denote the rational variety gotten
from $\PP^{d-1}$ by blowing up the points
$\ell_1,\ell_2,\ldots,\ell_n$. We write $L$ for the pullback of the
hyperplane class from $\PP^{d-1}$ to $X_G$ and $E_1,\ldots,E_n$ for
the exceptional divisors of the blow-up. The vector space
(\ref{sectionspace1}) can be rewritten as a space of sections:
\begin{equation}
 \label{sectionspace2}
  I_{(r,u)}^\perp \quad = \quad
 H^0\bigl(\,X_G, \,r L + (u_1{-}r) E_1 +  \cdots
+ (u_n{-}r) E_n\bigr),
\end{equation}
Taking the direct sum of these spaces for all $(r,u)$ in $\Z^{n+1}$,
which we identify with the Picard group of $X_G$, we obtain the
{\em Cox ring} of the blow-up:
$$ {\rm Cox}(X_G) \quad = \quad
\bigoplus_{(r,u) \in \Z^{n+1}} \!\! H^0\bigl(\,X_G, \,r  L +
(u_1{-}r) E_1 +  \cdots + (u_n{-}r) E_n\bigr). $$ We summarize our
discussion as follows:

\begin{cor} If $d \geq 3$ then the Cox-Nagata ring $R^G$ equals the Cox ring
of the variety $X_G$ which is gotten from $\PP^{d-1}$ by blowing up
the points $\ell_1,\ldots,\ell_n$.
\end{cor}

The Cox-Nagata ring $R^G$ has received a considerable amount of
attention in the recent literature in algebraic geometry. Some
relevant references are \cite{BP, CT, Der2, LV, Mukai, Vera, STV, TVV}.
These papers are primarily concerned with the case when
$\ell_1,\ldots,\ell_n $ are generic. In this case, the function
$\psi$ is invariant under the action of the Weyl group of the Dynkin
diagram $T_{2,d,n-d}$ with three legs of length $2, d$ and $n-d$.
The action of this Weyl group on $\Z^{n+1}$ is generated by
permutations of $ (u_1,\ldots,u_n)$ and the
 transformation $(r,u) \mapsto (r',u')$ where
\begin{eqnarray*}
& r' \,=\, \sum_{i=1}^d u_i - r \, , \,\,\,
 u_j' \,=\, u_j \,\,\,\hbox{for} \,\,\, j = 1,\ldots,d ,\\ &
u_j' \, = \, \sum_{i=1}^d u_i - 2r + u_j \,\,\,\hbox{for}\,\, j =
d+1,\ldots,n.
\end{eqnarray*}
While the invariance $\,\psi(r, u_1,\ldots,u_n) = \psi(r,u_{\pi(1)},
\ldots,u_{\pi(n)})\,$ for all permutations $\pi$ is entirely obvious
from the definition of $\psi$, the second invariance
\begin{equation}
\label{cremona}
 \psi(r,u_1,\ldots,u_n) \,\,\,= \,\,\, \psi(r',u_1',\ldots,u_n')
 \end{equation}
 is less obvious and due to Nagata \cite{Nag}.
 He proved (\ref{cremona})  by applying a Cremona transformation to the
points $\ell_1,\ldots,\ell_d$ in $\PP^{d-1}$. This induces an
isomorphism of the blow-up $X_G$ which replaces the divisor class
$\,r \cdot L + \sum_{j=1}^n (r-u_j) E_j\,$ by the divisor class
$\,r' \cdot L + \sum_{j=1}^n (r'-u'_j) E_j$. An alternative proof
using representation theory of $SL_2(K)$ was given by Kuttler and
Wallach \cite{KW}. The reader may find it instructive to verify
(\ref{cremona}) in Examples \ref{ex:eins}, \ref{ex:zwei} and
\ref{ex:drei}.

Our primary objective is to describe nice bases for $R^G$ which
express $\psi$ as the number of lattice points in a
$(d-1)$-dimensional polytope parametrized by $(r,u)$. The volume of
this polytope is the leading form of $\psi$. Algebraic geometers
know this as the {\em volume of the line bundle} $\,r \cdot L +
\sum_{i=1}^n (u_i{-}r) E_i\,$ on $X_G$. For instance, the area of
the polygon in Example \ref{ex:drei} would be an extension of the
formula in \cite[Example 3.5]{BKS} to blowing up $n=5$ points.

The {\em support semigroup} $\,\mathcal{S}^G\,$ of the Cox-Nagata
ring $R^G$ is the subsemigroup of $ \N^{n+1}$ which consists of all
multidegrees $(r,u)$ for which $\psi(r,u) > 0 $. The {\em support
cone} $\,\mathcal{C}^G \,$ is the cone in $\R^{n+1}$ generated by
the support semigroup.

\begin{remark} The support cone $\mathcal{C}^G$
is a full-dimensional cone in $\R^{n+1}$.
\end{remark}

If $\ell_1,\ell_2,\ldots,\ell_n$ are generic then the Weyl group
$T_{2,n-d,d}$ acts on the semigroup $\mathcal{S}^G$ and its cone
$\mathcal{C}^G$. It is known that $T_{2,n-d,d}$ is finite if and
only if (\ref{fracineq}) if and only if $R^G$ is finitely generated
if and only if its semigroup $\mathcal{S}^G$ is finitely generated
\cite{CT,  Mukai}. In this finite situation, the cone
$\mathcal{C}^G$ is the cone over an $n$-dimensional polytope
$\mathcal{P}^G$, which we call the {\em support polytope}. For $d=3$
and $n=6,7,8$ our Weyl group is $E_6,E_7,E_8$ respectively, and the
support polytopes are the beautiful {\em Gosset polytopes}
 \cite{gosset} associated with  these exceptional groups.
 However, in our view,
the support polytope and all the other concepts introduced in this
section are combinatorially interesting even if the $\ell_i$ are not
generic. Here is an example to illustrate this perspective.

\begin{example} \label{sixspecial}
We consider the six special points in the plane $\PP^2$ given by
$$ A \quad = \quad
\begin{bmatrix}
\,1 \,& \,0 \,& \,0\, & \phantom{-}1 & -1 & \phantom{-}0 \, \\
\,0 \,& \,1 \,& \,0 \,&            -1        & \phantom{-}0 & \phantom{-} 1 \,\\
\,0\, & \,0 \,& \,1 \,&  \phantom{-}0   &\phantom{-}1 & -1 \, \\
\end{bmatrix}.
$$
These are the intersection points of four general lines in $\PP^2$.
Here $d=3,n=6$, $\,G = {\rm kernel}(A)$, and $(\ell_1,\ldots,\ell_6)
= (z_1,z_2,z_3) \cdot A$. The Cox-Nagata ring equals
$$ R^G \,\,=\,\, K\bigl[x_1,x_2,x_3,x_4,x_5,x_6, L_{124}, L_{135}, L_{236}, L_{456},
M_{16}, M_{25}, M_{34}\bigr], $$ where the $L_{ijk}$ and $M_{ij}$
represent lines spanned by triples and pairs of points:
$$
\begin{matrix}
L_{124} & = & \underline{y_3  x_5  x_6} +x_3  y_5  x_6-x_3  x_5  y_6
\qquad \qquad \qquad &
{\rm in \ degree} \,\, (1001011)  ,\\
L_{135} & = & \underline{y_2  x_4  x_6} - x_2  y_4  x_6+x_2  x_4
y_6 \qquad \qquad \qquad  &
{\rm in \ degree} \,\, (1010101)  ,\\
L_{236} & = & \underline{y_1  x_4  x_5}+x_1  y_4  x_5-x_1  x_4  y_5
\qquad \qquad \qquad &
{\rm in \ degree} \,\, (1100110)  ,\\
L_{456} & = & \underline{y_1  x_2  x_3}+x_1  y_2  x_3+x_1  x_2  y_3
\qquad \qquad \qquad &
{\rm in \ degree} \,\, (1111000)  ,\\
 M_{16} & = & \underline{y_2  x_3  x_4 x_5} +x_2  y_3  x_4  x_5-x_2  x_3  y_4  x_5+x_2  x_3  x_4  y_5  &
{\rm in \ degree} \,\, (1011110)  ,\\
 M_{25} & = & \underline{y_1  x_3  x_4  x_6} +x_1  y_3  x_4  x_6+x_1  x_3  y_4  x_6-x_1  x_3  x_4  y_6 &
{\rm in \ degree} \,\, (1101101)  ,\\
 M_{34} & = & \underline{y_1  x_2  x_5  x_6} +x_1  y_2  x_5  x_6-x_1  x_2  y_5  x_6+x_1  x_2  x_5  y_6 &
{\rm in \ degree} \,\, (1110011) . \\
\end{matrix}
$$
The $7$-dimensional support semigroup $\mathcal{S}^G$ is spanned by
these seven vectors the unit vectors $e_i = {\rm degree}(x_i)$ for
$i=1,2,\ldots,6$. The $6$-dimensional support polytope
$\mathcal{P}^G$ has $13$ vertices. Its f-vector equals $(13, 69,
186, 260, 168, 38)$.
The seven underlined monomials together with $x_1,\ldots,x_6$ are
a sagbi basis of $R^G$, as defined in the next section.
The algorithm explained in Section~4 now computes an
Ehrhart-type formula for the Hilbert function $\psi$ of $R^G$.
 It outputs that $\psi(r,u)$ is the number of
integer vectors $(x,y) \in \N^2$ satisfying
\begin{eqnarray*}
&
x\,\leq {\rm min}(u_3,u_5,u_6), \,\,\,
   y\, \leq\, {\rm min}(u_2,u_4) ,\,\,\,
 2x+y \,\leq \, {\rm min}(u_3+u_6,u_5+u_6), \\ &
 3x+y\leq u_3+u_5+2u_6-r,\,\, r-u_2-u_4\leq x-y\leq u_3+u_5-r , \\ &
 \hbox{and} \,\,\quad
  r-u_1 \,\leq\, x+y\, \leq  \, {\rm min}(r,u_2 + u_5 -r,u_3+u_5+u_6-r).
\end{eqnarray*}
This piecewise quadratic function counts the sections of line bundles on a
singular surface
$X_G$ which is known 
classically as {\em Cayley's cubic surface}. \qed
\end{example}

\section{Sagbi bases}

In this section we introduce sagbi bases into the study of
Cox-Nagata rings, and we give a complete classification of such
bases when $n=d+1$. The basic idea is as follows. Let $K$ be the
field $\Q(t)$ of rational function in one variable $t$. As always in
tropical geometry \cite{speyer}, any field with a non-trivial non-archimedean
valuation would work as well, but to keep things as simple and
Gr\"obner-friendly as possible, we take $K= \Q(t)$ to be our field
of definition.

The {\em order} of a scalar $c(t) $ in $ K$ is the unique integer
$\omega$ such that $\,t^{-\omega} \cdot c(t)\,$ has neither a pole
nor a zero at $t=0$. Likewise, any polynomial $ f $ in $R = K[x,y]$
has a lowest order $\omega$ in the unknown $t$. We define the {\em
initial form} ${\rm in}(f)$ to be the coefficient of that lowest
power $t^\omega$ in $f$. In symbols,
$$ {\rm in}(f)  \,\,\, := \,\,\,
\bigl(t^{-\omega} \cdot f)|_{t=0} \quad \in \quad \Q[x,y]. $$ A
subset $\mathcal{F}$ of $K[x,y]$ is called {\em moneric} if ${\rm
in}(f)$ is a monomial for all $f \in \mathcal{F}$.

For any subalgebra $U$ of the polynomial ring $K[x,y]$, the {\em
initial algebra} ${\rm in}(U)$ is the subalgebra of $\Q[x,y]$
generated by the initial forms ${\rm in}(f)$ where $f$ runs over all
polynomials $f$ in $U$. Note that ${\rm in}(U)$ is a $\Q$-algebra
while $U$ is a $K$-algebra. Even if $U$ is a finitely generated
$K$-algebra, it will often happen that the $\Q$-algebra ${\rm
in}(U)$ is not finitely generated. A subset $\mathcal{F}$ of $U$ is
called a {\em sagbi basis} if $\mathcal{F}$ is {\em moneric} and the
initial algebra ${\rm in}(U)$ is generated as a $\Q$-algebra by the
monomials ${\rm in}(f)$ for $f \in \mathcal{F}$. The sagbi basis
$\mathcal{F}$ can be infinite even if $U$ is finitely generated. The
acronym {\em sagbi}
 was coined by Robbiano and Sweedler \cite{RS}. It stands for
 ``subalgebra analogue to Gr\"obner bases for ideals''.
 Our definition of sagbi bases is more general than the
definition usually given in the computer algebra literature
\cite{gobel, KM, MS}. There one uses a monomial order to define
initial monomials and the initial algebra, but this situation can be
modeled by introducing an extra variable $t$ as in \cite[\S
15.8]{Eis} to get to the situation described above. See also
\cite[\S 1]{CHV}.

\begin{remark} \label{whatweuse}
Throughout this paper we make frequent use of the basics
concerning the construction and properties
of sagbi bases which remain valid in our setting.
Such basics are:
{(a)} sagbi bases induce binomial initial ideals
  for partial term orders (as in
  \cite[Theorem 14.16]{MS}),
{(b)}   the lifting of all binomial syzygies implies the sagbi property
  (as in \cite[Proposition 1.3]{CHV}), 
{(c)} a multigraded algebra $R$ and its initial algebra ${\rm in}(R)$ share
Hilbert function. 
\end{remark}

\begin{example}
The example of a sagbi basis best known among computer algebraists
 is the set of elementary symmetric polynomials. Let $n=3$ and
 $\, \mathcal{F} \,\, = \,\, \bigl\{\, x_1 + t x_2 + t^2 x_3\,,\,\, x_1
x_2 + t x_1 x_3 + t^2 x_2 x_3 \,, \,\, x_1 x_2 x_3\, \bigr\}$.
The initial algebra of $\,U = K[\mathcal{F}]\,$
is $\, {\rm in}(U) =  \Q[ x_1, x_1 x_2, x_1
x_2 x_3] $.
If we wish to extend
to the invariants of the alternating group, then $U'$ is the $K$-algebra
generated by $\mathcal{F}$ and the discriminant $\,
(x_1 - t x_2) (x_1 - t^2 x_3) (x_2 - t x_3)$. Now, the initial
algebra $\,{\rm in}(U')\,$ is not finitely generated. See the
work of G\"obel \cite{gobel} for details.
 \qed
\end{example}

We seek to construct explicit sagbi bases of
Cox-Nagata subrings $R^G$  of $R = K[x,y]$. A first example was seen
in Example \ref{sixspecial} where the $13$ generators of $R^G$ form
a sagbi basis with the underlined initial monomials. The following
example shows that our sagbi bases do not generally
come from term orders.

\begin{example}
\label{FourPointsOnLine} Let $d = 2,n = 4$. Then $R^G =
K[x_1,x_2,x_3,x_4,E_1,E_2,E_3,E_4]$,
\vskip -.2cm
$$ \begin{matrix}
 \hbox{where} \qquad \qquad &
E_1 \,\,\, =\,\,  p_{23}   x_2 x_3 y_4 - p_{24}   x_2 y_3 x_4 +
p_{34}   y_2 x_3 x_4 ,\\ & E_2 \,\,\,= \,\,\,  p_{13}  x_1 x_3 y_4 -
p_{14}   x_1 y_3 x_4 + p_{34}   y_1 x_3 x_4 , \\ & E_3 \,\,\,=
\,\,\, p_{12}  x_1 x_2 y_4 - p_{14}   x_1 y_2 x_4 + p_{24}   y_1 x_2
x_4 , \\ & E_4 \,\,\, = \,\,\,  p_{12}  x_1 x_2 y_3 - p_{13}   x_1
y_2 x_3 + p_{23}   y_1 x_2 x_3.
\end{matrix}
$$
The coefficients $p_{ij}$ are any scalars in $K$ that satisfy
$\,p_{12} p_{34} - p_{13} p_{24}+  p_{14} p_{23} = 0\,$ and they
represent the Pl\"ucker coordinates of the subspace $G$ as in
(\ref{plucker}). The ring $R^G$ has Krull dimension $6$ and is
presented by the ideal of relations
\begin{eqnarray*} I^G \quad =  & \langle \,
p_{14} x_1 E_1 -  p_{24} x_2 E_2 + p_{34} x_3 E_3 ,\,
   p_{13}  x_1 E_1 - p_{23} x_2 E_2 + p_{34} x_4 E_4 ,  \\ & \,\,\,\,
     p_{12} x_2 E_2 - p_{13} x_3 E_3 + p_{14} x_4 E_4 ,\,
    p_{12} x_1 E_1 -   p_{23} x_3 E_3+ p_{24} x_4 E_4 \,\rangle.
\end{eqnarray*}
Among the $81=3^4$ choices of one term from each $E_i$, only $12$
are induced by a term order on $R$, and these are all equivalent
under permuting $\{1,2,3,4\}$. A representative is given by taking
the first term in each $E_i$ as listed above. These leading terms
specify a sagbi basis both in the classical sense and in our sense.
However, there is a another symmetry class  (also having $12$
elements) which specifies a sagbi basis in our sense but {\bf not}
via a term order in $(x,y)$. A representative is given by taking the
second term in each $E_i$. The corresponding initial toric algebras
of the
 Cox-Nagata ring $R^G$ are
 $$  \begin{matrix} {\rm in}(R^G) & = &
\Q[ \, x_1, \, x_2, \, x_3, \,  x_4, \,
 x_2 x_3 y_4,\,
 x_1 x_3 y_4, \,
 x_1 x_2 y_4,\,
 x_1 x_2 y_3\,]
 \\ & = &
\Q[x_1,{\ldots},x_4,E_1,{\ldots},E_4]/\langle x_1 E_1-x_2 E_2,x_2 E_2-x_3 E_3 \rangle \\
&&\hbox{e.g.~for}\,\,\, (p_{12},p_{13},p_{14},p_{23},p_{24},p_{34})
= (1,2,t,1,t,t),
\end{matrix}
$$
$$
  \begin{matrix} {\rm in}(R^G) & = &
  \Q[ \, x_1, \, x_2, \, x_3, \,  x_4, \,
x_2 y_3 x_4, x_1 y_3 x_4, x_1 y_2 x_4, x_1 y_2 x_3 \,]
\\ & = &
\Q[x_1,{\ldots},x_4,E_1,{\ldots},E_4]/\langle
x_1 E_1-x_2 E_2,x_3 E_3-x_4 E_4 \rangle \\
&& \!\!\!\! \hbox{e.g.~for} \, (p_{12},\ldots,p_{34}) =
(t^2-t^4, t-t^5, 1- t^6, t^2 - t^4, t- t^5, t^2{-}t^4).
\end{matrix}
$$
For both types we verify the sagbi property via \cite[Proposition
1.3]{CHV}, by observing that the binomial relations $x_i E_i - x_j
E_j$ lift to polynomials in $I^G$.
\qed
\end{example}

\begin{remark}
For $d \geq 3$, the Cox-Nagata ring $R^G$ serves as the  affine
coordinate ring  of the {\em universal torsor} \cite{Der2} over the
variety $X_G$. The geometric interpretation of sagbi bases
is that the universal torsor degenerates to a toric variety. By
restricting to the graded components of $R^G$, this induces
a simultaneous degeneration of each projective embedding of $X_G$ to
a projective toric variety. For instance, each lattice polygon in
Example \ref{ex:drei} represents a projective toric surface along
with a deformation to a del Pezzo surface. \qed
\end{remark}

We now present the classification of all sagbi bases for the
Cox-Nagata ring when $n = d+1$. 
Here we set $\,G = {\rm span}_K \{(\alpha_1,\ldots,\alpha_n)\} \simeq K^1\,$ in $\,K^n$.

\begin{thm} \label{grassthm}
Let $\mathcal{F}$ be the set consisting of  the $2 \times 2$-minors
of the matrix
$$ \begin{bmatrix}
\alpha_1 x_1 & \alpha_2 x_2 & \cdots & \alpha_n x_n \\
 y_1 & y_2 & \cdots & y_n
 \end{bmatrix}   $$
 and the variables $x_1,\ldots,x_n$.  Then  $\mathcal{F}$ is moneric if and only if
 the orders of the $\alpha_i$ are distinct. In this case
 $\mathcal{F}$ is a sagbi basis of the Cox-Nagata ring~$R^G$.
\end{thm}

\begin{proof}
After relabeling we may assume
$ \,{\rm ord}(\alpha_1)  >  {\rm ord}(\alpha_2) >  \cdots  >  {\rm ord}(\alpha_n)$.
The $2 \times 2$-minors have the initial monomials $x_i y_j$ for $1
\leq i < j \leq n$, and
\begin{equation}
\label{diagofminors} {\rm in}(\mathcal{F}) \,\, = \,\, \bigl\{
x_1,\ldots,x_n, \, y_1 x_2, y_1 x_3, \ldots, y_{n-1} x_n \,\bigr\}.
\end{equation}
That these monomials suffice to generate the initial algebra ${\rm
in}(K[\mathcal{F}])$ appears in
 \cite[\S 14.3]{MS}. Clearly, the elements of $\mathcal{F}$ are invariant under the
action of $G$ on $R$. The result that $R^G$ is generated by
$\mathcal{F}$ is classical \cite[Remark 3.9]{CT}.
\end{proof}

The Cox-Nagata ring $R^G$ is isomorphic to the algebra generated by
the $2 \times 2$-minors of a general $2 \times (n{+}1)$-matrix, that
is, the coordinate ring of the Grassmannian $Gr(2,n{+}1)$ of lines
in $\PP^n$. Theorem \ref{grassthm} represents
  the familiar sagbi degeneration of the Grassmannian
$Gr(2,n{+}1)$ to a toric variety.

We seek good formulas for evaluating the Hilbert function $\psi$ of
the Cox-Nagata ring $R^G$. For the case $n = d {+} 1$  discussed
here, we can utilize standard tools from
algebraic combinatorics. Recall (e.g.~from  \cite[\S 14.4]{MS}) that
a {\em  two-row Gelfand-Tsetlin pattern} is a non-negative integer
$2 \times n$-matrix $(\lambda_{ij})$ that satisfies  $\lambda_{2n} =
0$, $\lambda_{1,j+1} \geq \lambda_{2,j}$ and $\lambda_{i,j} \geq
\lambda_{i,j+1}$ for $i =1,2 $ and $j=1,\ldots,n-1$. Two-row
Gelfand-Tsetlin patterns are identified with monomials in our initial
algebra ${\rm in}(R^G)$ via the isomorphism of \cite[Theorem
14.23]{MS}. Here the monomials in (\ref{diagofminors}) correspond to
two-row partitions as in \cite[Example 14.21]{MS}. Using  this
isomorphism, our sagbi basis implies the following formula:

\begin{prop}\label{pr:d1}
For $n = d+1$, $\,\psi(r,u)$ is the number of two-rowed
Gelfand-Tsetlin patterns with $\lambda_{21} = r$ and $\lambda_{1j} +
\lambda_{2j} = u_j + \cdots + u_n$ for $j=1,\ldots,n$.
\end{prop}

This formula shows that $\psi$ is a piecewise polynomial function of
degree $n-2$ on the $(n+1)$-dimensional support cone
$\mathcal{C}^G$.
 The underlying $n$-dimensional support polytope
$\mathcal{P}^G$ is affinely isomorphic to the second hypersimplex
$\,\Delta(n{+}1,2) = {\rm conv} \{ e_i + e_j : 0 \leq i < j \leq n \}$.
 The $\Z^{n+1}$-grading specifies a linear map from a
Gelfand-Tsetlin polytope onto the second hypersimplex. Its fibers
are $(n-2)$-polytopes. They represent toric degenerations of all
projective embeddings of the blow-up of $\PP^{n-2}$ at $n$ points.

\begin{example}
In \cite[end of \S 4.1]{DI} the authors asked the question whether
\vskip -.2cm
$$
\varphi(j,\ldots,j)\,\,\, :=\,\,\, \sum_{r=0}^\infty\!\psi(r,j,\ldots,j) \,\, = \,\,
{\rm dim}_K (S/I_{(j,j,\ldots,j)})
$$
is a polynomial function in the one variable $j$.
We find a negative answer to this open question by examining
the case $d=2$ and $n=3$. By Proposition
\ref{pr:d1}, $\varphi(j,j,j)$ equals the Ehrhart quasi-polynomial of the
rational quadrangle
\begin{eqnarray*}
\{(\lambda_{11},\lambda_{12}, \lambda_{21},\lambda_{22})\in
\R_{\geq 0}^4&&:\,\,\, \lambda_{11}+\lambda_{21}=3,\,\,
\lambda_{12}+\lambda_{22}=2,\\
& & \lambda_{11}\geq \lambda_{12}\geq 1,\,\, \lambda_{21}\geq
\lambda_{22},\,\, \lambda_{12}\geq \lambda_{21},\,\,1\geq
\lambda_{22} \}.
\end{eqnarray*}
The number of lattice points in $j$ times this quadrangle equals
$$ \varphi(j,j,j) \,\, = \,\, \sum_{r \geq 0} \bigl( \, r+1 - 3 (r-j)_+ \bigr)_+ \,\,= \,\,
 \begin{cases}
(3j^2+ 6 j + 3)/4 & \hbox{ if $j$ is odd, }\\
(3j^2 + 6 j  + 4)/4 & \hbox{ if $j$ is even.}
\end{cases}
$$
This is not a polynomial in $j$ but it is a
 quasi-polynomial with period two. \qed
\end{example}

\section{Five points in the plane}

Let $G$ be a linear subspace of codimension $d$ in $K^n$. We call
the subspace $G$ {\em moneric} if the Cox-Nagata ring $R^G$ has a
minimal generating set $\mathcal{F}$ which is moneric, and we say
that $G$ is {\em sagbi} if $R^G$ has a minimal generating set
 $\mathcal{F}$ which is also a sagbi basis.
Two moneric subspaces $G$ and $G'$ are called equivalent if the
initial subalgebras ${\rm in}(R^G)$ and $ {\rm in}(R^{G'})$ are
identical and we call their common equivalence class a {\em moneric
class}. A moneric class can either be sagbi or not sagbi, depending
on whether the subspaces in that class are sagbi or not. In this
language, Theorem \ref{grassthm} says that for $n=d+1$, there is
precisely one moneric class up to permutations and this class is
sagbi.

We fix $n=d+2$ and we consider a generic linear subspace
\begin{equation}
\label{twoplane} G \,=\, {\rm rowspan}
\begin{bmatrix}
b_{11} & b_{12} & b_{13} & \cdots  & b_{1n} \\
b_{21} & b_{22} & b_{23} & \cdots  & b_{2n}
\end{bmatrix} \,\,\, \subset \,\, \, K^n.
\end{equation}
 Genericity means in particular  that the $\binom{n}{2}$  Pl\"ucker coordinates
\begin{equation}
\label{plucker}
 p_{ij} \,\,\, := \,\,\, b_{1i} b_{2j} - b_{1j} b_{2i} \qquad \qquad (1 \leq i < j \leq n)
 \end{equation}
are all non-zero. These Pl\"ucker coordinates already appeared in
Example \ref{FourPointsOnLine}, where the following result was
established: {\em There are precisely $24$ moneric classes of
$2$-dimensional linear subspaces of $K^4$. Modulo permutations of
$\{1,2,3,4\}$ the number of classes is two, and both of them are
sagbi classes.} Our main result in this section is the analogous
classification in the next case.

\begin{thm} \label{allbut60}
There are precisely $600$ moneric classes  of $3$-dimensional
generic linear subspaces of $K^5$. All but $60$ of these classes are sagbi. Modulo
permutations of the indices $\{1,\ldots,5\}$, the number of moneric classes is
seven and the number of sagbi classes is six.
\end{thm}

\begin{proof}
The Cox-Nagata ring $R^G$ is minimally generated by a distinguished
set $\mathcal{F}$ of $16$ polynomials \cite{BP, CT, STV}. First,
 the five variables $x_1,x_2,x_3,x_4,x_5$ represent
the exceptional divisors of the blow-up. Second, there are ten
generators corresponding to the lines passing through pairs of
points:
$$
L_{ijk} \quad := \quad p_{ij} \cdot x_i x_j y_k \,- \, p_{ik} \cdot
x_i y_j x_k \,+\, p_{jk} \cdot y_i x_j x_k \qquad \hbox{for} \,\,\,
1 \leq i<j<k \leq 5 . $$ And, finally, there is one generator for
the quadric through the five points:
$$
Q_{12345} \quad := \quad p_{12} p_{13} p_{23} p_{45} \cdot x_1 x_2
x_3 y_4 y_5 \,-\, p_{12} p_{14} p_{24} p_{35} \cdot x_1 x_2 y_3 x_4
y_5 \,\pm\, \cdots \cdots.
$$
The ten monomials $x_i x_j x_k y_l y_m$ in this expression
correspond to the ten splits $\{i,j,k\}\cup \{l,m\}$ of the index
set $\{1,2,3,4,5\}$. There are $3^{10} \times 10$ ways to chose one
monomial from each generator. To identify those choices that arise
from moneric subspaces $G$, we examine the {\em tropical Pl\"ucker
coordinates}
$$ d_{ij} \,\,:= \,\, - \,{\rm ord}(p_{ij}) \qquad \hbox{  for } \quad 1 \leq i < j \leq 5 .$$
After adding a positive constant, these quantities are non-negative,
and they represent the distances in a finite metric space. We order them
as follows:
$$  d \,\,= \,\, (d_{12},d_{13},d_{14},d_{15},
d_{23}, d_{24},d_{25},d_{34}, d_{35}, d_{45})  . $$ The metrics $d$
which arise from $2$-dimensional subspaces $G$ of $K^5$ are the
lattice points in
 the {\em tropical Grassmannian} ${\rm Trop}(Gr(2,5))$. This is a
$7$-dimensional fan in $\R^{10}$ whose combinatorial structure is
the Petersen graph \cite{speyer}. Each of its $15$ maximal cones
 is described as follows up to relabeling:
 \begin{equation}
 \label{Grass25cone}
   d_{ij} + d_{kl} = d_{ik} + d_{jl} \geq d_{il} + d_{jk}
\qquad \hbox{for} \,\, 1 \leq i < j< k < l \leq 5 .
\end{equation}
A computation reveals that the choice of leading terms for the
generating set $\mathcal{F}$ of $R^G$  divides the cone
(\ref{Grass25cone})
 into $160$ smaller convex cones, and thus we obtain  a
 finer fan structure with $2400$ maximal cones  on ${\rm Trop}(Gr(2,5))$.
 These $2400$ cones determine $600$ distinct monomial sets ${\rm in}(\mathcal{F})$.
  Modulo permuting the indices  $\{1,2,3,4,5\}$, the $600$ moneric classes
 come in  seven types:

\smallskip \noindent
{\em Type 1:} $d = (1, 2, 3, 4, 3, 4, 5, 1, 2, 3)$, ${\rm
in}(\mathcal{F}) = \{ x_1,x_2,x_3,x_4,x_5, y_1 x_2 x_3,
    y_1 x_2 x_4,$ $ y_1 x_2 x_5, x_1 y_3 x_4, x_1 y_3 x_5,
     x_1 y_4 x_5,       x_2 y_3 x_4,
     x_2 y_3 x_5, x_2 y_4 x_5, y_3 x_4 x_5, y_1 x_2 y_3 x_4 x_5 \} $,
     {\bf 120}.

\smallskip \noindent
{\em Type 2:} $d =  (5, 3, 5, 6, 4, 6, 7, 2, 3, 5)   $, ${\rm
in}(\mathcal{F})  = \{ x_1,x_2,x_3,x_4,x_5, x_1 x_2 y_3, y_1 x_2
x_4, $ $ y_1 x_2 x_5, x_1 y_3 x_4, x_1 y_3 x_5,
     x_1 y_4 x_5, x_2 y_3 x_4, x_2 y_3 x_5, x_2 y_4 x_5, y_3 x_4 x_5, y_1 x_2 y_3 x_4 x_5\}$,{\bf 120}.

\smallskip \noindent
{\em Type 3:} $d =    (2, 1, 3, 4, 3, 5, 6, 2, 3, 5)  $, ${\rm
in}(\mathcal{F})  = \{ x_1,x_2,x_3,x_4,x_5, y_1 x_2 x_3, y_1 x_2
x_4, $ $ y_1 x_2 x_5, x_1 y_3 x_4, x_1 y_3 x_5,
     y_1 x_4 x_5, x_2 y_3 x_4,  x_2 y_3 x_5, x_2 y_4 x_5, y_3 x_4 x_5, y_1 x_2 y_3 x_4 x_5\}$,{\bf 120}.

\smallskip \noindent
{\em Type 4:} $d =    (1, 1, 4, 4, 2, 5, 5, 3, 3, 6)  $, ${\rm
in}(\mathcal{F})  = \{ x_1,x_2,x_3,x_4,x_5, y_1 x_2 x_3, y_1 x_2
x_4, $ $ y_1 x_2 x_5, x_1 y_3 x_4, x_1 y_3 x_5,
     y_1 x_4 x_5, x_2 y_3 x_4, x_2 y_3 x_5, y_2 x_4 x_5, y_3 x_4 x_5, y_1 x_2 y_3 x_4 x_5\}$,
{\bf 60}.

\smallskip \noindent
{\em Type 5:} $d =     (1, 4, 5, 5, 5, 6, 6, 7, 7, 8) $, ${\rm
in}(\mathcal{F})  = \{ x_1,x_2,x_3,x_4,x_5, y_1 x_2 x_3, y_1 x_2
x_4,$ $ y_1 x_2 x_5, y_1 x_3 x_4, y_1 x_3 x_5,
     y_1 x_4 x_5, y_2 x_3 x_4, y_2 x_3 x_5, y_2 x_4 x_5, y_3 x_4 x_5, y_1 y_3 x_2 x_4 x_5 \} $,
{\bf 60}.

\smallskip \noindent
{\em Type 6:} $d =   (1, 1, 2, 3, 2, 3, 4, 1, 2, 1)   $, ${\rm
in}(\mathcal{F})  = \{ x_1,x_2,x_3,x_4,x_5, y_1 x_2 x_3, y_1 x_2
x_4, $ $ y_1 x_2 x_5, x_1 y_3 x_4, x_1 y_3 x_5,
     x_1 y_4 x_5, x_2 y_3 x_4, x_2 y_3 x_5, x_2 y_4 x_5, x_3 y_4 x_5, y_1 x_2 x_3 y_4 x_5\}$,
{\bf 60}.

\smallskip \noindent
{\em Type 7:} $d =   (1, 2, 3, 3, 3, 4, 4, 5, 5, 6)  $, ${\rm
in}(\mathcal{F})  = \{x_1,x_2,x_3,x_4,x_5, y_1 x_2 x_3, y_1 x_2 x_4,
$ $y_1 x_2 x_5, y_1 x_3 x_4, y_1 x_3 x_5,
     y_1 x_4 x_5, y_2 x_3 x_4,  y_2 x_3 x_5, y_2 x_4 x_5, y_3 x_4 x_5, y_1 y_2 x_3 x_4 x_5 \} $,
     {\bf 60}.

\medskip
We now check for each of the seven moneric types whether it is
sagbi. This is done as follows. We first compute the {\em toric
ideal} of algebraic relations among the $16$ monomial generators.
For Types 1,2,3,4,5 and 6, this toric ideal is minimally generated
by $20$ quadratic binomials, while Type 7 requires $21$ binomial
generators. This immediately implies that Type 7 is not sagbi,
because one the $21$ binomials does not lift to a relation on $R^G$,
thus violating the sagbi criterion in \cite[Proposition 1.3]{CHV}.
For Types 1--6 the $20$ binomials are distributed over ten different
multidegrees. For instance, for Type 6 the minimal generators form a
Gr\"obner basis (with leading terms listed first):
$$
\begin{matrix}
\hbox{degree} && \qquad \hbox{two generators} \!\!\!\!\!\! & \!\!\!\!\!\!\!\! \hbox{in that degree} \quad \\
(2,1,1,1,1,2) &&  L_{125} L_{345} - x_5 Q_{12345}, & L_{145} L_{235} - L_{135} L_{245} \\
(2,1,1,1,2,1) &&  L_{124} L_{345} - x_4 Q_{12345}, &  L_{145} L_{234} - L_{134} L_{245} \\
(2,1,1,2,1,1) && L_{123} L_{345} - x_3 Q_{12345} ,&   L_{135} L_{234} - L_{134} L_{235} \\
(2,1,2,1,1,1) && L_{123} L_{245} - x_2 Q_{12345} ,&  L_{125} L_{234} - L_{124} L_{235} \\
(2,2,1,1,1,1) && L_{123} L_{145} - x_1 Q_{12345}, &  L_{125} L_{134} - L_{124} L_{135} \\
(2,0,1,1,1,1) && x_3 L_{245} - x_2 L_{345} ,& x_5 L_{234} - x_4 L_{235} \\
(2,1,0,1,1,1) && x_3 L_{145} - x_1 L_{345}, & x_5 L_{134} - x_4 L_{135} \\
(2,1,1,0,1,1) && x_2 L_{145} - x_1 L_{245} ,&  x_5 L_{124} - x_4 L_{125} \\
(2,1,1,1,0,1) && x_2 L_{135} - x_1 L_{235} ,& x_5 L_{123} - x_3 L_{125} \\
(2,1,1,1,1,0) && x_2 L_{134} - x_1 L_{234}, & x_4 L_{123} - x_3 L_{124}. \\
\end{matrix}
$$
Each of these binomials lifts to a relation in the Cox-Nagata ring.
For instance, the relation $\, x_2 L_{134} - x_1 L_{234}\,$  among
the
 underlined leading monomials of
$\, L_{134}\, =\, p_{14} \cdot \underline{x_1 y_3 x_4} +\, \cdots
\,\,$
 and $\,L_{234} \,=\, p_{24} \cdot \underline{x_2 y_3 x_4} + \,\cdots \,$  lifts to the relation
$\,  p_{24}  x_2  L_{134} - p_{14} x_1 L_{234} - p_{34}  x_3
L_{124} \,$ in the prime ideal $I^G$ of relations on $R^G$.  Indeed,
the twenty generators of $I^G$ and their degrees are well-known, and
the existence of these quadrics verifies the sagbi property for
types 1--6.
\end{proof}

\begin{cor}
The Cox-Nagata ring $R^G$ is normal, Gorenstein and Koszul.
\end{cor}

\begin{proof}
This result is known by work of Popov
 \cite{popov}. Our approach offers a combinatorial certificate.
The Gr\"obner basis for the toric ideal of Type 6 is squarefree and
quadratic. We conclude that the toric algebra ${\rm in}(R^G)$ is
normal and Koszul, and hence so is its flat deformation $R^G$. 
 The Gorenstein property of $R^G$ holds because
the toric algebras ${\rm in}(R^G) $ of Types 1--5 are Gorenstein
with the same $\Z^6$-graded Betti numbers as $R^G$.
\end{proof}

The Cox-Nagata ring $R^G$ and its initial algebras ${\rm in}(R^G)$
share the same Hilbert function $\psi$. Each of the six sagbi types
specifies a formula for the piecewise quadratic function $\psi$. For
instance, the formula in Example \ref{ex:drei}  comes from Type 6.
We explain how this works.  The variables $y_2$ and $y_5$ are absent
from  ${\rm in}(\mathcal{F})$. Hence every monomial in the
subalgebra ${\rm in}(R^G)$ of $\Q[x,y]$
  has the form  $\, x^a y^b =
x_1^{a_1} x_2^{a_2} x_3^{a_3} x_4^{a_4} x_5^{a_5} y_1^{b_1}
y_3^{b_3} y_4^{b_4} \,$ and corresponds to a lattice point $(a,b) =
(a_1,a_2,a_3,a_4,a_5,b_1,b_3,b_4)$ in $\Z^8$. We now perform a
change of variables $\, (a,b) \mapsto (r,u,x,y)\,$ in the lattice
$\Z^8$ as follows:
$$
\begin{matrix}
r = b_1 + b_3 + b_4, &  x = b_1, &  y = b_4 , &
u_1 = a_1 + b_1 , \\
u_2 = a_2 + b_2 , & u_3 = a_3 + b_3 , & u_4 = a_4 + b_4 , & u_5 =
a_5 + b_5.
\end{matrix}
$$
Let $\Gamma$  denote the convex cone in $\,\R^{8}$ generated by the
$16$ vectors $(r,u,x,y)$ corresponding to  the $16$ monomials in
${\rm in}(\mathcal{F})$. Then $\,\Gamma \,\cap \, \Z^{8}\,$ is the
normal semigroup whose semigroup algebra equals ${\rm in}(R^G)$. Our
degree function is now the projection
  $(r,u,x,y) \mapsto (r,u)$ onto the first six coordinates, and
     the Hilbert function of  our Type 6 semigroup $\Gamma \cap \Z^8$
  for this grading equals
  \begin{equation}
  \label{nicepsiformula}
   \psi(r,u) \,\,\,= \,\,\,
\# \bigl\{ (x,y) \in \Z^2 \,:\, \,(r,u,x,y) \in \Gamma  \bigr\}  .
\end{equation}
What is left to do is to compute the cone $\Gamma$. The software
{\tt polymake} reveals that
 the f-vector of $\Gamma$ equals
 $(16 ,80 ,180, 216 ,148, 58, 12)$,
 and it generates the twelve facet-defining inequalities
 of $\Gamma$. These are precisely the
 twelve inequalities, such as $y-x \leq u_4 + u_5 - r$,
 which are listed in Example \ref{ex:drei}.

We had chosen Type 6 for the formula in Example \ref{ex:drei}
because it gives the fewest linear inequalities (only $12$). For the
other types, the $8$-dimensional cone $\Gamma$ has $\geq 14$ facets.
More precisely, $f(\Gamma) = (16, 84, 200 , 253 , 180 ,71 ,14)$  for
Type 1, and $f(\Gamma) = (16, 87, 221 ,301, 229, 94 ,18)$ for Types
2, 3, 4 and 5.

The $540$ cones $\Gamma$ all share the same image under the linear
map
$${\rm degree} \, :\, \R^8 \,\rightarrow \R^6 \, ,  \,\, (r,u,x,y) \mapsto (r,u). $$
That image is the $6$-dimensional support cone $\,\mathcal{C}^G =
{\rm degree}(\Gamma)$, with  f-vector $(16,80,160,120,26)$.
The underlying support polytope $\,\mathcal{P}^G\,$ is the
{\em demicube}
$$ \mathcal{P}^G \,\, = \,\, {\rm conv}
\bigl((10000), \ldots, (00001), (11100), (11010), \ldots, (00111),
(11111) \bigr\}.$$ Each of our $540$ sagbi bases specifies an
$7$-dimensional polytope $\Pi$ with a distinguished projection
$\,\Pi \rightarrow \mathcal{P}^G\,$ onto the $5$-dimensional demicube. Its
fibers over the interior of  $\mathcal{P}^G$ are $2$-dimensional
polygons, one for each ample line bundle on the del Pezzo
surface $X_G$. The corresponding projective embedding of $X_G$
degenerates to the toric surface associated with that polygon.

\section{The cubic surface}

Let $G$ be generic linear subspace of dimension $3$ in $K^6$ and let
$p_{ijk}$ denote its dual Pl\"ucker coordinates, i.e., $p_{ijk}$ is
the $3 \times 3$-minor with column indices $(i.j,k)$ of a $3 \times
6$-matrix $A$ whose kernel equals $G$. The Cox-Nagata ring $R^G$ has
$27$ minimal generators. They correspond to the $27$ lines in the
cubic surface gotten from $\PP^2$ by blowing up the points in the
columns of $A$. There are six generators $E_i = x_i$ representing
the exceptional divisors, $15$ generators $F_{ij}$ representing the
lines through pairs of points, and  six generators $G_i$
representing the quadrics through any five of the points. We can
express these as polynomials in $x$ and $y$ whose coefficients are
expressions in the $p_{ijk}$. For instance, the generator for the
line through points $1$ and $2$  equals
$$ F_{12} \,\,\, = \,\,\,
p_{123} \cdot y_3 x_4 x_5 x_6 + p_{124} \cdot x_3 y_4 x_5 x_6 +
p_{125} \cdot x_3 x_4 y_5 x_6 + p_{126} \cdot x_3 x_4 x_5 y_6 ,$$
and the generator for the quadric through the points $1,2,3,4,5$
equals
\begin{eqnarray*}
G_6 \,\,\,\, = &     p_{123} p_{124} p_{125} p_{345}  \cdot  y_1
y_2   x_3   x_4   x_5   x_6^2
  \, \, -\,\, p_{132} p_{135} p_{134} p_{245}    \cdot   y_1   y_3   x_2   x_4   x_5   x_6^2 \\ & \!\!\!\!\!\!\!\!
\,  \, + \,\, p_{142} p_{143} p_{145} p_{235}   \cdot   y_1   y_4
x_2   x_3   x_5   x_6^2
 \, \, - \,\, p_{152} p_{153} p_{154} p_{234} \cdot     y_1   y_5   x_2   x_3   x_4   x_6^2 \\ &  \!\!\!\!\!\!\!\!
\, \,  + \,\, p_{231} p_{234} p_{235} p_{145}    \cdot  y_2   y_3
x_1   x_4   x_5   x_6^2
 \, \, - \,\, p_{241} p_{243} p_{245} p_{135}  \cdot    y_2   y_4   x_1   x_3   x_5   x_6^2 \\ &  \!\!\!\!\!\!\!\!
\, \,  + \,\, p_{251} p_{253} p_{254} p_{134}  \cdot    y_2   y_5
x_1   x_3   x_4   x_6^2 \, \,  + \, \, p_{341} p_{342} p_{345}
p_{125}  \cdot    y_3   y_4   x_1   x_2   x_5   x_6^2 \\ &
\!\!\!\!\!\!\!\! \,  \, - \,\, p_{351} p_{352} p_{354} p_{124}
\cdot   y_3   y_5   x_1   x_2   x_4   x_6^2
\, \,  + \,\, p_{451} p_{452} p_{453} p_{123}  \cdot     y_4   y_5   x_1   x_2   x_3   x_6^2 \\
 &
\,\,  + \,\, (p_{124} p_{235} p_{136} p_{145}-p_{123} p_{245}
p_{146} p_{135})
  \cdot y_1   y_6   x_2   x_3   x_4   x_5   x_6 \\ &
+ \,\, (p_{124} p_{135} p_{236} p_{245}-p_{123} p_{145} p_{246}
p_{235}) \cdot  y_2   y_6   x_1   x_3   x_4   x_5   x_6 \\ & + \,\,
(p_{134} p_{125} p_{236} p_{345}+p_{123} p_{145} p_{346} p_{235})
\cdot  y_3   y_6   x_1   x_2   x_4   x_5   x_6 \\ & + \,\, (p_{124}
p_{135} p_{346} p_{245}-p_{134} p_{125} p_{246} p_{345}) \cdot  y_4
y_6   x_1   x_2   x_3   x_5   x_6 \\ & +\,\, (p_{125} p_{134}
p_{356} p_{245}-p_{135} p_{124} p_{256} p_{345}) \cdot  y_5   y_6
x_1   x_2   x_3   x_4   x_6 \\ & + \,\,\,(p_{124} p_{135} p_{236}
p_{456}-p_{123} p_{145} p_{246} p_{356})  \cdot
  y_6^2   x_1   x_2   x_3   x_4   x_5 . \,\,\,
\end{eqnarray*}
Note that $\,{\rm deg}(F_{12}) = (1,0,0,1,1,1,1)\,$ and $\,{\rm
deg}(G_6) =  (2,1,1,1,1,1,2)$.

We shall now construct a toric model for the cubic surface $X_G$ as
follows.

\begin{thm}\label{th:cox6sagbi}
There exists a $3$-dimensional sagbi subspace $G$ of $K^6$ whose
 toric ideal ${\rm in}(I^G)$ is generated by quadrics
and has a squarefree Gr\"obner basis.
\end{thm}

\begin{proof}
Let $G$ denote the kernel of the matrix
$$ A \quad = \quad
\begin{bmatrix}
\, t^{1} & t^{11} & t^{11} & t^{13} & t^{7} & t^{7}  \, \\
\, t^6      & 1  & t^{13} & t^{10} & t^{15} & t^{15} \, \\
\, t^{6} & t^{5} & 1 & t^{15} & t^{5} & t \,
\end{bmatrix}.
$$
This subspace induces the following initial monomials for the $27$
generators:
$$
\begin{matrix}
\,E_1 \,& \,E_2 \,& \,E_3 \,& \,E_4 \,& \,E_5 \,& \,E_6 \,& F_{12} & F_{13} & F_{14} & F_{15} \\
x_1 & x_2 & x_3 & x_4 & x_5 & x_6  & y_3x_4x_5x_6&
 y_2x_4x_5x_6 &
x_2y_3x_5x_6 & y_2x_3x_4x_6
\end{matrix}
$$
$$
\begin{matrix}
 F_{16} & F_{23} & F_{24} & F_{25} & F_{26} & F_{34}  & F_{35} \\
  y_2x_3x_4x_5&
 y_1x_4x_5x_6 &
 x_1y_3x_5x_6 &
y_1x_3x_4x_6 & y_1x_3x_4x_5 &
 y_1x_2x_5x_6 &
x_1y_2x_4x_6
      \end{matrix}
 $$

$$
\begin{matrix}
  F_{36} & F_{45} & F_{46} & F_{56} & G_6 & G_5 \\
  x_1y_2x_4x_5  &y_1x_2x_3x_6&
y_1x_2x_3x_5& x_1y_2x_3x_4 & x_1y_2y_3x_4x_5x_6^2&
x_1y_2y_3x_4x_5^2x_6
\end{matrix}
$$
$$
\begin{matrix}
 G_4 & G_3 & G_2 & G_1 \\
  x_1y_2y_3x_4^2x_5x_6&
x_1y_2y_3x_3x_4x_5x_6&  x_1x_2y_2y_3x_4x_5x_6& x_1^2y_2y_3x_4x_5x_6.
 \end{matrix}
$$
The toric ideal of relations among these $27$ monomials is minimally
generated by $81$ quadrics.
 These generators occur as triples in $27$ distinct
degrees. For example,  in degree $\, (1,1,1,1,0,1,1)\,$ we find the
three generators $\,
E_{3}F_{34}-E_{6}F_{46},\,E_{1}F_{14}-E_{2}F_{24},\,E_{3}F_{34}-E_{5}F_{45}$,
in  degree  $\, (3,2,1,2,2,2,2)  \,$ we find the three minimal
generators $\,F_{12}G_{1}-F_{24}G_{4},\,F_{23}G_{3}-F_{25}G_{5},\,
F_{25}G_{5}-F_{26}G_{6}\,$ etc.

Each of these $81$ binomial relations among our $27$ monomials
lifts to a quadratic relation in the presentation ideal $I^G$ of the Cox-Nagata ring $R^G$.
The lifting property can be checked either directly, by writing down an quadratic
polynomial in $I^G$ whose
initial form is the given binomial, or indirectly, by computing the values of $\psi$
on each of the $27$ observed degrees:
$$ \psi(1,1,1,1,0,1,1) = 2, \,\,\,
\psi (3,2,1,2,2,2,2) = 2, \, \ldots \ldots. $$ Of course, only one such
computation suffices if we use the fact that all $27$ degrees are in
a single orbit under the Cremona action (\ref{cremona}) of the Weyl
group $E_6$. Using  \cite[Proposition 1.3]{CHV} we conclude that
${\rm in}(R^G)$ is generated by the $81$ monomials listed above, and
that the toric ideal ${\rm in}(I^G)$ is generated by the $81$
quadratic binomials $\, E_{3}F_{34}-E_{6}F_{46},\ldots,
F_{25}G_{5}-F_{26}G_{6}, \ldots$.

A computation shows that
the reduced Gr\"obner basis of the toric ideal
${\rm in}(I^G) $ with respect to the
reverse lexicographic order has squarefree initial monomials.
This property ensures that the toric algebra ${\rm in}(R^G)$ is normal,
and we can expect a nice polyhedral formula for its Hilbert function  $\,\psi $.
\end{proof}

Just like in Section 4, our sagbi basis automatically implies
a polyhedral formula for $\psi$, namely, we need to
compute the facet inequalities of the cone $\Gamma$ underlying
the normal toric  algebra ${\rm in}(R^G)$. Using {\tt polymake}, we find
$$ f(\Gamma) \,\,= \,\, (27, 216, 747, 1287, 1191, 603, 162, 21). $$
The $21$ facet inequalities translate into a formula for $\psi$
if we perform the following change of variables
$(a,b)\mapsto (r,u,x,y)$ in the lattice $\Z^9$:
 \begin{eqnarray*}
 r=b_1+b_2+b_3,\,\,\,\,\,\,\, x=b_2,\,\,\,\,\, y=b_3,\,\,\,\,\,
 u_1=a_1+b_1,\\
 u_2=a_2+b_2,u_3=a_3+b_3,u_4=a_4,u_5=a_5,u_6=a_6.
 \end{eqnarray*}
 After this change of coordinates,  $\Gamma$
 is the convex cone in $\R^9$ generated by the $27$
 vectors $(r,u,x,y)$ corresponding to the 27 monomial generators
 $x^ay^b=x_1^{a_1}x_2^{a_2}x_3^{a_3}x_4^{a_4}x_5^{a_5}x_6^{a_6}y_1^{b_1}y_2^{b_2}y_3^{b_3}$
 of ${\rm in}(R^G)$. Then formula
(\ref{nicepsiformula}) holds, and we find:

\begin{cor} \label{twentyone}
If $d=3$, $n=6$ and the linear forms $\ell_1,\ell_2,\ldots,\ell_6$ are generic then
 $\psi(r,u_1,u_2,\ldots,u_6)$ equals the number of lattice points $(x,y)$ satisfying
 \begin{eqnarray*}
 & \!\!\!\!\!
 \max(0,2r{-}u_2{-}u_3{-}u_4,2r{-}u_2{-}u_3{-}u_6)\leq
\! x \! \leq
\min(u_1,u_5,u_1{+}u_4{+}u_5{+}u_6{-}2r), \\
  & \!\!\!\!
  \max(0,2r {-} u_1 {-} u_4 {-} u_6,2r {-} u_4 {-} u_5 {-} u_6)\leq \! y \! \leq
 \min(u_2,u_3, u_1{+}u_2{+}u_3{+}u_5 {-} 2r ),\\
  &\!\!\!\!\!\! \max(r - u_4,r - u_6,3r - u_2 - u_3 - u_4 - u_6)\leq x+y\leq \\ &
  \qquad \qquad \qquad \qquad
   \min(r,u_1+u_3+u_5-r,u_1+u_2+u_5-r),\\
  & r-u_2-u_3\leq x-y,\,\,\,
 2x+y \leq u_1+u_5,  \,\,\, \hbox{and} \,\,\,\,
2r-u_4-u_6\leq x+2y .
\end{eqnarray*}
\end{cor}

\begin{remark}
Using the formula above, we can rapidly compute the dimension
of the space of sections (\ref{sectionspace1}) of any line bundle
on the del Pezzo surface $X_G$ of degree three. In particular, we can check
that $\psi(3,2,2,2,2,2,2) = 4$, which corresponds
to the familiar embedding
of $X_G$ as a cubic surface in $\PP^3$.
\end{remark}

We conjecture that the formula of Corollary \ref{twentyone} is optimal,
in the sense that the number $21$ of linear inequalities is minimal.
Equivalently, the number of facets of the $9$-dimensional cone  $\Gamma$
arising from any  three-dimensional sagbi subspace $G$ of $K^6$ is at least $21$.
A proof of this conjecture might require the complete classification of the
equivalence classes of  subspaces, as was done for $n = 5$ in
Section 4. We suggest this as a research problem:

\begin{problem} \label{prob:naruki}
Determine all equivalence classes of $3$-dimensional sagbi subspaces of $K^6$,
i.e.~extend the classification of Theorem \ref{allbut60} from $K^5$ to $K^6$.
\end{problem}

In such a classification, the role of the tropical Grassmannian
${\rm Trop}(Gr_{2,5})$ would be played by a suitable
tropical model of the moduli space of cubic surfaces.
We expect this model to be a variant of the fan constructed by Hacking, Keel
and Tevelev \cite{HKT}. The Naruki coordinates on
that moduli space correspond to the coefficients of
the generators of $R^G$. An example of such a coordinate
is the coefficient $\,p_{124} p_{135} p_{236}
p_{456}-p_{123} p_{145} p_{246} p_{356} \,$ appearing
in the generator $G_6$ listed prior to Theorem  \ref{th:cox6sagbi}.
The task of Problem \ref{prob:naruki} is to divide that variant of
the Hacking-Keel-Tevelev fan into many small cones,
one for each system of choices of leading terms for the
$27$ generators of $R^G$.

\section{Del Pezzo surfaces of degree one and two}

The computational result of the previous section
extends to $n=7$ and $n=8$:

\begin{thm}\label{th:cox8sagbi}
Let $4 \leq n \leq 8$. Then there
exists a generic $3$-dimensional sagbi subspace $G$ of $K^n$ whose
 toric ideal ${\rm in}(I^G)$ is generated by quadrics.
\end{thm}

Since ${\rm in}(I^G)$ is a flat degeneration of $I^G$, Theorem \ref{th:cox8sagbi} implies in
particular that the presentation ideal $I^G$ of the Cox-Nagata ring $R^G$ is
generated by quadrics. This furnishes a
computational proof of the following result  which
was obtained independently also by
Testa, V\'arilly-Alvarado and Velasco \cite{TVV}.

\begin{cor}
The presentation ideal of the Cox ring of
a del Pezzo surface gotten by blowing up at most $8$
general points in $\PP^2$ is generated by quadrics.
\end{cor}

This quadratic generation result had been conjectured by
Batyrev and Popov in \cite{BP},
and it was proved for $n \leq 7$
in the subsequent papers \cite{Der1, LV, Vera}.

\begin{proof}[Proof of Theorem 6.1]
The case $n=4$ is covered by Theorem \ref{grassthm}
and the cases $n=5$ and $n=6$ are dealt with in
Theorems \ref{allbut60} and \ref{th:cox6sagbi}
respectively. In what follows we present
choices of sagbi subspaces $G$ for $n=7$ and $n=8$.

First consider the case $n=7$ where $X_G$ is
 the Del Pezzo surfaces of degree  two.
 The Cox-Nagata ring $R^G$ has 56
minimal generators.  There are seven generators
$E_i=x_i$ for the exceptional divisors, $21$ generators
$F_{ij}$ representing the lines through pairs $\{\ell_i,\ell_j\}$ of points, $21$
generators $G_{ij}$ representing the quadrics through any five of
the points, missing $\{\ell_i,\ell_j\}$, and seven generators $C_i$
representing the cubics through all seven points, where $\ell_i$ is a
double point.

Let $G$ be the $4$-dimensional subspace of $K^7$ which is the
kernel of
\begin{equation}
\label{MatrixFor7}
 A \quad = \quad
\begin{bmatrix}
t^3&t^{10}&1&t^6&t^{17}&t^{12}&t^{11}\\
t^{18}&t^{15}&t^8&t^4&t^6&t^7&t\\
t^{10}&t^{16}&t^2&1&t^6&t^4&t^9
\end{bmatrix}.
\end{equation}
This subspace induces the following initial monomials for the $56$
generators:
$$
\begin{matrix}
E_1 & E_2 & E_3 & E_4 & E_5 & E_6 &\ E_7& F_{12} & F_{13} & F_{14}  \\
x_1 & x_2 & x_3 & x_4 & x_5 & x_6  & x_7& x_3x_5x_6x_7y_4 &
x_2x_4x_5x_6y_7 & x_2x_3x_5x_6y_7
\end{matrix}
$$
$$
\begin{matrix}
 F_{15}&  F_{16} & F_{17} & F_{23} & F_{24} & F_{25}   \\
x_2x_3x_6x_7y_4& x_2x_3x_4x_5y_7& x_2x_3x_5x_6y_4& x_1x_4x_5x_6y_7&
x_1x_3x_5x_6y_7 & x_1x_3x_6x_7y_4
\end{matrix}
 $$
$$
\begin{matrix}
F_{26}&  F_{27}& F_{34} & F_{35} & F_{36} & F_{37}    \\
 x_1x_3x_4x_5y_7&x_1x_3x_5x_6y_4& x_1x_2x_5x_6y_7&
x_1x_2x_6x_7y_4 & x_1x_2x_4x_5y_7 & x_1x_2x_5x_6y_4
\end{matrix}
$$
$$
\begin{matrix}
F_{45}&  F_{46}&F_{47} & F_{56} & F_{57} & F_{67} \\
x_1x_2x_6x_7y_3&x_1x_2x_5x_7y_3&x_1x_2x_5x_6y_3  & x_1x_2x_4x_7y_3 &
x_1x_2x_4x_6y_3&
 x_1x_2x_4x_5y_3
 \end{matrix}
$$
$$
\begin{matrix}
G_{67}& G_{57} & G_{56} & G_{47}  \\
x_1x_2x_3x_5x_6^2x_7y_4y_7&x_1x_2x_3x_4x_5^2x_6y_7^2 &
x_1x_2x_3x_5^2x_6^2y_4y_7&x_1x_2x_3x_4x_5x_6x_7y_4y_7
\end{matrix}
$$
$$
\begin{matrix}
G_{46}&G_{45}&G_{37}&G_{36}\\
x_1x_2x_3x_5x_6^2x_7y_4^2&x_1x_2x_3x_4x_5^2x_6y_4y_7&x_1x_2x_3x_5x_6x_7^2y_3y_4&
x_1x_2x_3x_5x_6^2x_7y_3y_4
\end{matrix}
$$
$$
\begin{matrix}
G_{35}&G_{34}&G_{27}&G_{26}\\
x_1x_2x_3x_5^2x_6x_7y_3y_4
&x_1x_2x_3x_4x_5x_6x_7y_3y_4&x_1x_2^2x_4x_5x_6x_7y_3y_7&
 x_1x_2^2x_4x_5x_6^2y_3y_7
\end{matrix}
$$
$$
\begin{matrix}
G_{25}&G_{24}&G_{23}&G_{17}\\
x_1x_2^2x_4x_5^2x_6y_3y_7 &x_1x_2^2x_4^2x_5x_6y_3y_7 &
x_1x_2^2x_3x_4x_5x_6y_3y_7&
 x_1^2x_2x_4x_5x_6x_7y_3y_7
\end{matrix}
$$
$$
\begin{matrix}
G_{16}&G_{15}&G_{14}&G_{13}\\
x_1^2x_2x_4x_5x_6^2y_3y_7 & x_1^2x_2x_4x_5^2x_6y_3y_7 &
x_1^2x_2x_4^2x_5x_6y_3y_7& x_1^2x_2x_3x_4x_5x_6y_3y_7
\end{matrix}
$$
$$
\begin{matrix}
G_{12}&C_{1}&C_{2}&C_{3}\\
x_1^2x_2^2x_4x_5x_6y_3y_7 & \!\! x_1x_2^2x_3x_4x_5^2x_6^2x_7y_3y_4y_7 & \!\!
x_1^2x_2x_3x_4x_5^2x_6^2x_7y_3y_4y_7& \!\!
x_1^2x_2^2x_4^2x_5^2x_6^2y_3y_7^2
\end{matrix}
$$
$$
\begin{matrix}
C_{4}&C_{5}&C_{6}&C_{7}\\ \!\!
x_1^2x_2^2x_3x_4x_5^2x_6^2y_3y_7^2 & \!\!\!
x_1^2x_2^2x_3x_4x_5x_6^2x_7y_3y_4y_7 & \!\!\!
x_1^2x_2^2x_3x_4^2x_5^2x_6y_3y_7^2& \!\!\!
x_1^2x_2^2x_3x_4x_5^2x_6^2y_3y_4y_7.
\end{matrix}
$$
The toric ideal of relations among these $56$ monomials has $529$
minimal generators. All generators are quadratic binomials. This
computation is performed most efficiently using the software {\tt 4ti2}
due to Hemmecke {\it et al.} \cite{4ti2}. The $56$
monomials involve $10$ variables. We input these data as a $10
\times 56$-matrix into {\tt 4ti2} and apply the command {\tt
markov}. The resulting {\em Markov basis} consists of $529$ vectors,
each representing a quadratic binomial in $56$ unknowns.

To check that the subspace $G$ is sagbi, it remains to be verified
that each of the $529$ quadratic binomials lifts to a relation in $I^G$.
This is done in a self-contained manner by computing $\psi(r,u)$ -- directly
from the definition -- for each degree $(r,u)$ that occurs among the
$529$ binomials in the output from {\tt 4ti2}.
To verify the lifting condition of \cite[Proposition 1.3]{CHV}, it
suffices to check that $\psi(r,u)$ plus the number of binomials of
degree $(r,u)$ equals the number of all monomials of degree $(r,u)$
in the quantities $E_i, F_{ij}, G_{ij}$ and $C_i$.

Of course, this last computational step would be unnecessary if we
allow ourselves to apply to prior results from the literature. Indeed,
it was already shown by Derenthal \cite[\S 4]{Der1}
that $I^G$ contains precisely $529$ linearly independent
quadrics, so we just need to compare
the degrees output by {\tt 4ti2} with the degrees in Derenthal's generators,
and the sagbi property follows. However, we wish to emphasize that our
approach is independent of any prior work, as
we can  verify the sagbi property of $G$ by computing a few values of $\psi$.
We conclude that ${\rm in}(I^G)$ is generated by $529$ quadrics
and hence so is $I^G$.

\smallskip

We now come to the hardest case, $n=8$, where $X_G$ is the
del Pezzo surface of degree one. Derenthal \cite[Lemma 15]{Der1} computed
that the ideal $I^G$ contains $17399$ linearly independent quadrics,
and he showed that these quadrics generate $I^G$ up to radical.
Testa, V\'arilly-Alvarado and Velasco \cite{TVV} applied the methods of
\cite{LV} to this situation, and they succeeded in proving that these
$17399$ quadrics do indeed generate the prime ideal $I^G$.
Our sagbi approach gives an alternative proof which is computational
and elementary.
The only prior knowledge we are using is that $R^G$ has
$242$ minimal generators.

 Let $G$ be the $5$-dimensional subspace of $K^8$
which is the kernel of

\begin{equation}
\label{MatrixFor8}
 A \quad = \quad
\begin{bmatrix}
t^{6}&t^{10}&t^{3}&t^{10}&t^{1}&t^{4}&t^{10}&t^{2}\\
t^{10}&t^{3}&t^{8}&t^{6}&t^{8}&t^{1}&t^{8}&t^{8}\\
t^{4}&t^{8}&t^{7}&t^{7}&t^{8}&t^{5}&t&t^{9}
\end{bmatrix}.
\end{equation}
The initial form of each minimal generator of $R^G$
is found to be a monomial, so the subspace $G $ is moneric.
Consider this list of $242$ monomials in
$x_1,\ldots,x_8,y_1,\ldots,y_8$.
Actually, of the eight $y$-variables only three appear among
these monomials, so we are facing a list of $242$ monomials
in $11$ variables, and our task is to compute the toric ideal
of algebraic relations among these monomials.
We write the list of monomials as a $11 \times 242$-matrix of non-negative
integers, we input that matrix into the software {\tt 4ti2}, and we apply
the command {\tt markov} to compute minimal generators of the toric ideal.

After two days or so, the computation terminates. The output
is an integer matrix with $17399$ rows and $242$ columns. Each
row represents a binomial in our toric ideal, and we check
that all $17399$ binomials are quadrics.
Using the same technique  as  in the proof for $n=7$, namely computing a few values of $\psi$,
we verify that all these $17399$ quadratic binomials lift to polynomials in $I^G$.
This proves that both ${\rm in}(I^G)$ and $I^G$ are generated by quadrics.
\end{proof}

\begin{remark}
All input and output files for the software {\tt 4ti2} used in this proof are posted
at the website {\tt http://lsec.cc.ac.cn/$\sim$xuzq/cox.html}.
\end{remark}

Our proof of the Batyrev-Popov conjecture for del Pezzo surfaces
of degree one is a fairly automatic process, albeit
computational intensive. Given standard computer algebra tools
to carry out the verification of the sagbi property, the proof only amounts to
 exhibiting the matrices (\ref{MatrixFor7}) and (\ref{MatrixFor8}).
Naturally, it would be extremely interesting to explore the moduli space
of possibilities for such matrices,  along the lines suggested in Problem \ref{prob:naruki}.

\smallskip

Our motivation for this study was to find a formula for
the counting function $\psi$. The sagbi matrices
(\ref{MatrixFor7}) and (\ref{MatrixFor8}) specify explicit
Ehrhart-type formulas for $\psi$, just like in
Example \ref{ex:drei} and Corollary \ref{twentyone}.
The shape of that formula appears in (\ref{nicepsiformula}),
and it is made completely explicit by listing the facets of the
convex polyhedral cone $\Gamma$ associated with the
toric algebra ${\rm in}(R^G)$. For $n=7$ the cone $\Gamma$
 is $10$-dimensional and has $56$ rays.  For
 $n=8$ the cone $\Gamma$ is $11$-dimensional and has $240$ rays.
An even higher-dimensional example is discussed in
Example \ref{ex:ana}. See also
Corollary \ref{decoration} for  $n=d+2$.

When the Cox-Nagata ring $R^G$ is not finitely generated,
(e.g.~when the inequality (\ref{fracineq}) does not hold), 
the convex cone $\Gamma$
for ${\rm in}(R^G)$ still exists but it will no longer
be polyhedral. The equation (\ref{nicepsiformula}) remains
valid, but the question how to make this formula explicit
and useful requires further study.

\section{Phylogenetic algebraic geometry}

In this section we fix $n=d+2$ and we assume that $G$ is the
rowspace of a $2 \times n$-matrix $(b_{kl})$ as in (\ref{twoplane})
whose Pl\"ucker coordinates $p_{ij} = b_{1i} b_{2j} - b_{1j} b_{2i}$
are all non-zero.
 For any  subset
$\{ i_0 < i_1 < \ldots < i_{2k}\}$ of
 $\{1,2,\ldots,n\}$ having odd cardinality
we define $\,Q_{i_0 i_1 \cdots i_{2k}} \,$ to be the determinant of
the matrix
\begin{equation}
\label{CTmatrix}
\begin{bmatrix}
b_{1 i_0}^k x_{i_0} & b_{1 i_1}^k x_{i_1} & b_{1 i_2}^k x_{i_2} &
b_{1 i_3}^k x_{i_3} & \cdots &
b_{1 i_{2k}}^k x_{i_{2k}} \\
b_{1 i_0}^{k-1} y_{i_0} & b_{1 i_1}^{k-1} y_{i_1} & b_{1 i_2}^{k-1}
y_{i_2} & b_{1 i_3}^{k-1} y_{i_3} & \cdots &
b_{1 i_{2k}}^{k-1} y_{i_{2k}} \\
b_{1 i_0}^{k-1} b_{2 i_0} x_{i_0} & b_{1 i_1}^{k-1} b_{2 i_1}
x_{i_1} & b_{1 i_2}^{k-1} b_{2 i_2} x_{i_2} & b_{1 i_3}^{k-1} b_{2
i_3} x_{i_3} & \cdots &
b_{1 i_{2k}}^{k-1} b_{2 i_{2k}} x_{i_{2k}} \\
b_{1 i_0}^{k-2} b_{2 i_0} y_{i_0} & b_{1 i_1}^{k-2} b_{2 i_1}
y_{i_1} & b_{1 i_2}^{k-2} b_{2 i_2} y_{i_2} & b_{1 i_3}^{k-2} b_{2
i_3} y_{i_3} & \cdots &
b_{1 i_{2k}}^{k-2} b_{2 i_{2k}} y_{i_{2k}} \\
\vdots & \vdots & \vdots & \vdots & \ddots & \vdots \\
b_{2 i_0}^k x_{i_0} & b_{2 i_1}^k x_{i_1} & b_{2 i_2}^k x_{i_2} &
b_{2 i_3}^k x_{i_3} & \cdots &
b_{2 i_{2k}}^k x_{i_{2k}} \\
\end{bmatrix} \! .
\end{equation}
This is a homogeneous polynomial with
$$ {\rm deg}(Q_{i_0 i_1 \cdots i_{2k}}) \,\,\, = \,\,\,
 k e_0 + e_{i_0} + e_{i_1}  + \cdots + e_{i_{2k}}. $$
Moreover, $Q_{i_0 i_1 \cdots i_{2k}}$ is invariant under the Nagata
action. Indeed, the action  by any vector in $G$ adds to each
$y$-row in  the matrix (\ref{CTmatrix}) a linear combination of the
two adjacent $x$-rows. This leaves the determinant unchanged.
 We can write the coefficients of this polynomial
 as products of Pl\"ucker coordinates:
 \begin{equation}
 \label{expandedQ}
 Q_{i_0 i_1 \cdots i_{2k}}\,\, = \,\,
 \sum \pm \, (\!\! \prod_{i,j \in \{0,\ldots,k\}} \!\!\!\!\! p_{a_i a_j} \,)\cdot
( \!\! \! \prod_{r,s \in \{1,\ldots,k\}} \!\!\!\!\! p_{b_r,b_s}\,)
\cdot
 x_{a_0} x_{a_1} \cdots x_{a_k} y_{b_1} \cdots y_{b_k}
 \end{equation}
 where the sum is over all partitions
  $$ \{i_0,i_1,\ldots,i_{2k}\} \,\,= \,\, \{a_0,a_1,\ldots,a_k\} \, \cup \,
 \{b_1,\ldots,b_k\}.$$

 \begin{thm} {\rm (Castravet-Tevelev \cite{CT})} \label{CTTheorem}
The Cox-Nagata ring $R^G$ is minimally generated by the $2^{n-1}$
invariants $Q_{i_0 i_1 \cdots i_{2k}}$ where $1 \leq i_0 < \cdots <
i_{2k} \leq n$.
\end{thm}

This result is stated in \cite[Theorem 1.1]{CT} for the case when
one row of the matrix $(b_{ij})$ consists of ones. The general case
easily follows because the del Pezzo variety $X_G$ and its Cox ring
remain unchanged if the blown-up points in $\PP^{n-3}$ undergo a
projective transformation. It is important to note, however, that
our sagbi analysis will not work if the vector of ones lies in $G$.

Another key ingredient for this section is the celebrated {\em Verlinde formula}
which lies at the interface of algebraic geometry and mathematical physics.
We refer to equation (12.6) of Mukai's book \cite[\S 12]{Mukbook}.
The following interpretation of that formula and its proof were suggested
to us by Jenia Tevelev.

\begin{thm} {\rm (Verlinde formula)}
\label{verlinde}
For $n = d+2$ and $G$ generic we have
$$
\psi(\,dl\,,\,2l,2l,\ldots,2l) \quad = \quad
\frac{1}{2l+1} \sum_{j=0}^{2l} (-1)^{dj}
\bigl({\rm sin} \frac{2j+1}{4l+2} \pi \bigr)^{-d}.
$$
Here $l$ can be a half-integer if $d$ is even
but must be an integer if $d$ is odd.
\end{thm}

\begin{proof}
The right hand side is the number of sections
of multiples of the canonical line bundle on the
moduli space $\mathcal{N}_{0,n}$ of rank two
stable quasiparabolic vector bundles on $\PP^1$
with $n=d+2$ marked points. See \cite[\S 12.5]{Mukbook}.
A result due to Stefan Bauer \cite{Bau} states that
$\mathcal{N}_{0,n}$ and the blow-up
$X_G$ of $\PP^{n-3}$ at $n$ points are related by a sequence of flops.
See  \cite[\S 2]{Mukai}  and after
\cite[Theorem 12.56]{Mukbook}.
This implies that $\mathcal{N}_{0,n}$
and the blow-up $X_G$ have the same Picard group,
their Cox rings are isomorphic, and their canonical classes agree.

The anticanonical class on $X_G$ equals
$\, -K \, = \, dH-(d-2)(E_1+ \cdots +E_{d+2}) $.
This is a primitive element in the Picard group 
when $d$ is odd but divisible by $2$ when $d$ is even \cite[Remark 12.54]{Mukbook}.
The right hand sum above equals
$$ {\rm dim}_K H^0 \bigl( \mathcal{N}_{0,d+2}, \mathcal{O}(-lK) \bigr)
\quad = \quad {\rm dim}_K H^0 \bigl( X_G , \mathcal{O}(-lK) \bigr). $$
Our equation  (\ref{sectionspace1})
implies that this dimension
equals $\psi(r,u_1,\ldots,u_n)$ where
$\,r = dl,\,u_i-r = l(2{-}d)\,$ and $l$
is allowed to be a half-integer if $d$ is even.
\end{proof}

Hilbert functions of ideals generated by powers of $n$ general linear
forms in $n-2$ unknowns were studied by D'Cruz and Iarobbino in
\cite{DI}. The previous two theorems establish both parts of their
main conjecture on page 77 of \cite{DI}.

\begin{cor}  \label{DIconjodd}
 If $n = 2k+1$ is odd then $\psi(k,1,1,\ldots,1) = 1$.
\end{cor}

\begin{proof}
If $n=2k+1$ then  $\,Q_{1 2 \cdots n}\,$ is the unique
minimal generator of the Cox-Nagata ring $R^G$ in degree $(k,1,1,\ldots,1)$. 
This follows from Theorem \ref{CTTheorem},
due to Castravet and Tevelev.
Hence the
space $\,R^G_{(k,1,1,\ldots,1)}\,$ is one-dimensional.
\end{proof}

\begin{cor} \label{DIconjeven}
 If $n=2k+2$ is even then $\psi(k,1,\ldots,1) = 2^k$.
 \end{cor}
 
 \begin{proof}
 For $d=n-2 = 2k$ and $l=1/2$ the
 trigonometric sum in the Verlinde formula
 (Theorem \ref{verlinde}) simplifies to $2^k$.
 See also \cite[page 483, line 4]{Mukbook}.
 We are grateful to Jenia Tevelev for suggesting this derivation to us.
  \end{proof}
 
 We are now prepared to embark towards
 the punchline of this section: 
 {\em  Sagbi bases connect the Verlinde formula to
 phylogenetic algebraic geometry}.
Let $K[q]$ and $\Q[q] $ denote the polynomial rings over $K$ and
$\Q$ in $2^{n-1}$ variables $q_{i_0 i_1 \cdots i_{2k}}$, one for
each odd subset of $\{1,\ldots,n\}$. We seek to compute the
presentation  ideal $I^G \subset K[q]$ of the Cox-Nagata ring $R^G$.
Our approach to this problem is to identify subspaces $G$ that are
 sagbi. This allows us to study $I^G$ by way of its toric
initial ideal ${\rm in}(I^G) $ in $ \Q[q]$. We begin
 by introducing a class of toric ideals in $\Q[q]$. These
 represent statistical models in \cite{BW, SS}.

Let $T$ be a trivalent phylogenetic tree which leaves labeled by
$\{1,\ldots,n\}$. Each interior edge of $T$ corresponds to a {\em
split}, by which we mean an unordered pair $\{A,B\}$ such that $A
\cup B = \{1,\ldots,n\}$, $A \cap B = \emptyset$, $|A| \geq 2$ and
$|B| \geq 2$. The set of all $n-3$ splits of $T$ is denoted ${\rm
splits}(T)$. By \cite[Theorem 2.35]{ASCB}, the combinatorial type of
the tree $T$ is uniquely specified by the set ${\rm splits}(T)$.

For each $\{A,B\} \in {\rm splits}(T)$ we introduce two matrices
whose entries are variables in $\Q[q]$. They are denoted ${\bf
M}_{A,B}$ and ${\bf M}_{B,A}$. Their formats are $2^{|A|-1} \times
2^{|B|-1}$ and  $2^{|B|-1} \times 2^{|A|-1}$ respectively. For the
matrix ${\bf M}_{A,B}$, the rows are indexed by even subsets
$\sigma$ of $A$, the columns are indexed by odd subsets $\tau$ of
$B$, and the entry in row $\sigma$ and column $\tau$ is the variable
$q_{\sigma \cup \tau}$. Similarly, the entries of $M_{B,A}$ are the
variables $q_{\sigma \cup \tau}$ for $\sigma \subseteq B$ even and
$\tau \subseteq A$ odd. We write ${\bf I}_T$ for the ideal in
$\Q[q]$ which is generated by the $2 \times 2$-minors of all $2n-6$
matrices ${\bf M}_{A,B}$ and ${\bf M}_{B,A}$ where $\{A,B\}$ runs
over ${\rm splits}(T)$. Known results in phylogenetic algebraic
geometry \cite{BW} imply that ${\bf I}_T$ is a prime ideal:

\begin{thm} \label{BWtheorem}
The ideals ${\bf I}_T$ are toric and they all have the same Hilbert
function with respect to the $\Z^{n+1}$-grading $\, {\rm deg}(q_{i_0
i_1 \cdots i_{2k}}) \, = \,
 k e_0 + e_{i_0}  + \cdots + e_{i_{2k}}$.
 \end{thm}

\begin{proof}
Let $\Q[q']$ be the polynomial ring whose variables $q'_{j_1 j_2
\ldots j_{2k}}$ are indexed by the even subsets
$\{j_1,j_2,\ldots,j_{2k}\}$ of $\{1,2,\ldots,n\}$. We declare the
leaf $n$ to be the root of the tree $T$ and we identify $\Q[q]$ with
$\Q[q']$ by mapping $\, q_\sigma \mapsto q'_{\sigma \backslash
\{n\}}\,$ if $n \in \sigma\,$ and $\, q_\sigma \mapsto q'_{\sigma
\cup \{n\}}\,$ if $n \not\in \sigma$. The image of ${\bf I}_T$ in
$\Q[q']$ under this identification coincides with the prime ideal of
the binary symmetric model \cite{BW} (called {\em Jukes-Cantor
model} in \cite{ASCB, SS}) on the phylogenetic tree $T$. Indeed, it was
shown in \cite[\S 6.2]{SS} that the $2 \times 2$-minors of the above
matrices form a Gr\"obner basis for ${\bf I}_T$ with respect to a
suitable term order. Buczy\'nska and Wi\'sniewski
  \cite[Theorem 3.26]{BW} proved that all toric ideals ${\bf I}_T$
  for the various $T$ have
the same Hilbert function in the standard $\Z$-grading. However,
their proof works verbatim for our finer $\Z^{n+1}$-grading as well.
\end{proof}

The result by Buczy\'nska and Wi\'sniewski \cite{BW} that the Hilbert
function of ${\bf I}_T$ is invariant under any choice of trivalent tree
$T$ is remarkable because there are as many as $(2n-5)!!$ distinct
trees $T$. In \cite{BW} the question is left open whether the toric
ideals ${\bf I}_T$ all lie on the same irreducible component of the
corresponding multigraded Hilbert scheme, and, if yes, what is the general point
on that component. We here answer this question, by constructing
sagbi deformations of the toric varieties in \cite{BW} to
the projective variety with coordinate ring $R^G$.
Here it is essential that we use the definition of sagbi bases
given in Section 3. Sagbi bases for term orders in $K[x,y]$
will not work.

Let $\mathcal{F}$ be the set of the $2^{n-1}$ minimal generators
$Q_{i_0 i_1 \cdots i_{2k}}$ in Theorem \ref{CTTheorem}. Suppose that
$G$ is a moneric subspace of codimension $2$ in $K^n$. This means
that the $2^{n-1}$ elements ${\rm in}(Q_{i_0 i_1 \cdots i_{2k}})$ in
the set
 ${\rm in}(\mathcal{F})$ are all monomials.
 Let ${\bf J}_{G}$ denote the kernel
of the ring map $\pi:\Q[q] \rightarrow \Q[x,y]$ which takes the
variable $q_{i_0 i_1 \cdots i_{2k}}$ to the monomial ${\rm
in}(Q_{i_0 i_1 \cdots i_{2k}})$. In other words, ${\bf J}_{G}
\subset \Q[q]$ is the toric ideal of algebraic relations among the
initial monomials of $\mathcal{F}$. Suppose that $T$ is a trivalent
tree with leaves $1,2,\ldots,n$ and ${\bf I}_T$ its toric ideal as
above. We say that the moneric subspace $G$ {\em realizes} the tree
$T$ if $\,{\bf I}_T \, = \, {\bf J}_G$.

\begin{example} \label{caterpillar}
Let $G$ be the row space of the $2 \times n$-matrix
$$
[b_{ij}] \quad  = \quad
\begin{bmatrix}
     1 & t & t^2 & \cdots & t^{n-3} & t^{n-2} & t^{n-1} \\
\,t^{n-1} & t^{n-2} & t^{n-3} & \cdots & t^2 & t & 1
\end{bmatrix}.
$$
Then $G$ is moneric since, and for each odd subset $\{i_0 < \cdots <
i_{2k}\}$, we have
\begin{equation}
\label{diagmono}
 {\rm in} (Q_{i_0 i_1 \cdots i_{2k}}) \,\, = \,\,
x_{i_0} y_{i_1} x_{i_2} y_{i_3} x_{i_4} \, \cdots \, y_{i_{2k-1}}
x_{i_{2k}} .
\end{equation}
This is seen either from the matrix (\ref{CTmatrix}), whose
diagonal entries multiply to (\ref{diagmono}), or from the
expansion (\ref{expandedQ}).  Let $T$ be the {\em caterpillar tree} whose
splits~are
$$ {\rm Splits}(T) \,\, = \,\, \bigl\{
\{A,B\} \,: \, {\rm max}(A) < {\rm min}(B) \bigr\}. $$ We claim that
the subspace $G$ realizes the tree $T$. For each split $\{A,B\}$ of
$T$ we consider the images of the matrices ${\bf M}_{A,B}$ and ${\bf
M}_{B,A}$ under the map $\pi$. The matrix $\pi({\bf M}_{A,B})$
equals the product of the column vector labeled by even sets $\sigma
= \{\sigma_1 {<} \cdots {<} \sigma_{2i}\} \subseteq A$ and entries
$x_{\sigma_1} y_{\sigma_2} x_{\sigma_3} \cdots y_{\sigma_{2i}}$,
with the row vector labeled by odd sets $\tau = \{\tau_0 {<} \cdots
{<} \tau_{2i}\} \subseteq A$ and entries $x_{\tau_0} y_{\tau_1}
x_{\tau_2} \cdots y_{\tau_{2i}}$. This uses the assumption
$\sigma_{2i} < \tau_{0}$. The matrix $\pi({\bf M}_{B,A})$ is a
similar product with the roles of $B$ and $A$ reversed. Hence the
matrices ${\bf M}_{A,B}$ and ${\bf M}_{B,A}$ have rank one modulo
${\bf J}_G$, and this implies
 $\,{\bf J}_G \subseteq {\bf I}_T$.
Since both ideals are prime of the same Krull dimension, namely
$2n-2$, it follows that
 $\,{\bf J}_G = {\bf I}_T$. \qed
\end{example}

Our next lemma says that the caterpillar tree is not alone:

\begin{lem} \label{lemreal}
For any trivalent tree $T$ there is a subspace $G$ which
realizes~$T$.
\end{lem}

\begin{proof}
We proceed by induction on $n$. For $n \leq 5$ every trivalent tree
is a caterpillar tree, so we are done by Example \ref{caterpillar}.
Let $n \geq 6$ and fix any split of $T$. We cut the tree along the
split into two smaller trees $T'$ and $T''$, each having a leaf in
the place of that split. By induction, the trees $T'$ and $T''$ can
be realized by subspaces $G'$ and $G''$. Let $d'$ and $d''$ be the
corresponding metric spaces, given by negated orders of the
Pl\"ucker coordinates of $G'$ and $G''$. We build a new tree metric
by joining the tree metrics $d'$ and $d''$ along the split, but in
such a way that the length of the split edge is {\bf much larger}
than any of the previous edge lengths. The effect of this choice is
that any initial monomial ${\rm in}(Q_{i_0 i_1 \cdots i_{2k}})$ is
the product of previous initial monomials coming from $G'$ and $G''$
with a variable $x_{\rm split}$ removed where needed. Given the
metric $d$ we determine the initial monomials from the expansion
(\ref{expandedQ}). Our choice ensures that each matrix
$\pi(M_{A,B})$ or $\pi(M_{B,A})$ is a product of a column vector
times a row vector. As in Example \ref{caterpillar} we conclude
$\,{\bf J}_G \subseteq {\bf I}_T$ and, as both ideals are prime of
Krull dimension $2n{-}2$, it follows that
 $\,{\bf J}_G = {\bf I}_T$.
\end{proof}

Our next two lemmas establish the relationship to the
Verlinde formula.

\begin{lem} \label{lemAA}
 Let $u$ be a vector in
$\{0,1,2\}^n$ which has $i$ entries $0$, has $j$ entries $2$ and has
  $2k+2$ entries $1$.
 Then
  $\psi(j+k,u) = 2^k$.
\end{lem}

\begin{proof}
If $n=2k+2$ then this is the content of Corollary \ref{DIconjeven}.
Next consider the case when all entries in $u$ are $0$ or $1$. After
relabeling we may assume $u = (1,1,\ldots,1,0,\ldots,0)$, and we can
take the linear forms corresponding to the $0$'s to be variables,
say, $\, \ell_{2k+3} = z_{n-2k+1} ,\ldots, \ell_n = z_{n-2}$. Then
we have
$$ \,I_u \,\, = \,\,
\left\langle \ell_1^2, \ell_2^2, \ldots, \ell_{2k+2}^2,  z_{2k+1}
,\ldots, z_{n-2} \right\rangle, $$ and
$\psi(k,1,1,\ldots,1,0,\ldots,0) = \psi(k,1,1,\ldots,1) = 2^k$ is
clear from the definition of $\psi$. For the general case we
induction on $j$, and we apply the fact that $\psi$ is invariant
under the action of the  Weyl group $ D_n$ by  Cremona
transformation. Using this action on our degree $(j+k,u)$, we can
replace each entry $0$ in $u$ by an entry $2$ while incrementing the
first coordinate $j+k$ by one.
\end{proof}

\begin{lem} \label{lemBB}
Let $u$ be a vector in $\{0,1,2\}^n$ with $i$ entries $0$, $j$
entries $2$ and
  $2k+2$ entries $1$.
Then ${\rm dim}(\Q[q]_{(j+k,u)}) = 2^{2k}$ and ${\rm
dim}((\Q[q]/{\bf I}_T)_{(j+k,u)}) =  2^k$.
\end{lem}

\begin{proof}
The monomials of degree $(j+k,u)$ in $\Q[q]$ are products $q_\sigma
q_\tau$ where $|\sigma \cup \tau| = j+k$, $|\sigma \cap \tau| = j$,
$u_s = 2$ for $s \in \sigma \cap \tau$, and $u_t = 0$ for $t \not\in
\sigma \cup \tau$. The set $\,(\sigma \backslash \tau) \cup (\tau
\backslash \sigma)\,$ has $2k+2$ elements, and the above products
$q_\sigma q_\tau$ are in bijection with partitions of this set into
two odd subsets. There are $2^{2k}$ such partitions, and this
implies the first assertion $\,{\rm dim}(\Q[q]_{(j+k.u)}) =
2^{2k}$.

For the second assertion we may assume that $T$ is the caterpillar
tree, in light of Theorem  \ref{BWtheorem}, so we have ${\bf I}_T =
{\bf J}_G$ as in Example \ref{caterpillar}. Let us first consider
the case $i=j=0$. A monomial in the algebra generated by
(\ref{diagmono}) has degree $(k,u) = (k,1,1,\ldots,1)$ if and only
if it can be factored in the form
\begin{equation}
\label{factoredmono} x_1
\cdot (x_2 y_3 \,\,\hbox{or} \,\,y_2 x_3) \cdot (x_4 y_5
\,\,\hbox{or} \, \,y_4 x_5) \, \cdots \cdot \,(x_{2k} y_{2k+1}
\,\,\hbox{or} \,\,  y_{2k} x_{2k+1} ) \cdot x_{2k+2}.
\end{equation}
The number of distinct such products equals $2^k$ as required. The
general case is obtained by removing both variables $x_s$ and $y_s$
whenever $u_s = 0$, and by including both variables $x_t$ and $y_t$
in the above product when $u_t = 0$.
\end{proof}

We are now prepared to state and prove our main result in this section.

\begin{thm}
Every trivalent phylogenetic tree $T$ is realized by subspace $G$ which is sagbi, and 
the counting function $\psi$ equals the
common $\Z^{d+3}$-graded Hilbert function of the toric algebras
$\Q[q]/{\bf I}_T$ associated with  the trees~$T$.
\end{thm}

\begin{proof}
Using  Lemma \ref{lemreal} we find a subspace $G$ which realizes the
given tree $T$. The ideal of algebraic relations among the initial
monomials ${\rm in}(\mathcal{F})$ equals the toric ideal ${\bf
I}_T$. The generators of ${\bf I}_T$ are the $ 2 \times 2$-minors of
the matrices ${\bf M}_{A,B}$ and ${\bf M}_{B,A}$ for $\{A,B\} \in
{\rm Splits}(T)$. Lemma \ref{lemBB} describes the number of linearly
independent generators of ${\bf I}_T$ in each multidegree $(k,u)$.
Lemmas \ref{lemAA} and \ref{lemBB} tell  us that the Cox-Nagata ring
$R^G = K[\mathcal{F}]$ has the same number of relations in each such
degree $(k,u)$. This means that the inclusion $\,{\rm
in}(K[\mathcal{F}]_{(k.u)}) \subseteq \Q[{\rm
in}(\mathcal{F})]_{(k,u)}\,$ is an equality for each syzygy degree
$(k,u)$. This means that each binomial relation on ${\rm
in}(\mathcal{F})$ lifts to a relation on $\mathcal{F}$. Using
\cite[Proposition 1.3]{CHV}, we conclude that $\mathcal{F}$ is a
sagbi basis of $R^G$, which means that $G$ is sagbi. The statement
about the Hilbert function follows.
\end{proof}

\begin{cor}
The ideal $I^G$ of the Cox-Nagata ring
is generated by quadrics.
\end{cor}

\begin{proof}
The initial toric ideal 
${\rm in}(I^G)  = {\bf J}_G = {\bf I}_T$ is
generated by quadrics.
\end{proof}

We next derive an explicit piecewise polynomial formula
for the function $\psi : \N^{n+1} \rightarrow \N$. Consider any trivalent tree $T$ with $n$ leaves and let
$\rho$ be a positive integer. We define a {\em $T$-decoration} of
order $\rho$ to be an assignment of nonnegative integer weights to
the edges of $T$ such that the half-weight of every interior node is
an integer bounded above by  $\rho$ and bounded below by each
adjacent edge weight. Here the  {\em half-weight} of an interior
node of $T$ is defined as half the sum of the weights of the three
adjacent edges.

\begin{cor} \label{decoration}
Fix any trivalent phylogenetic tree $T$ with $n$ leaves. Then
$\psi(r,u_1,\ldots,u_n)$  equals the number
    of T-decorations of order $\rho =
    u_1 + u_2 + \cdots + u_n - 2 r \,$
    whose pendant edges have weights $u_1,u_2,....,u_{n-1},\rho-u_n$.
    \end{cor}

 \begin{proof}[Sketch of Proof]
 This uses the realization
  of the Jukes-Cantor ideal ${\bf I}_T$ in the polynomial  $\Q[q']$
 whose variables $q'_{j_1 j_2 \ldots j_{2k}}$
are indexed by the even subsets $\{j_1,j_2,\ldots,j_{2k}\}$. The
parametrization discussed in \cite{BW, SS} maps the variable
$q'_{j_1 j_2 \ldots j_{2k}}$ to the unique collection of edges in
$T$ which connect the leaves $j_1, j_2, \ldots, j_{2k}$ pairwise by
$k$ edge-disjoint paths on $T$. The subsemigroup of $\N^{{\rm
edges}(T)}$ generated by these sets of edges is saturated, and the
linear inequalities describing the corresponding cone are precisely
the inequalities in terms of half-weights of the interior nodes in
our definition of $T$-decorations.
 \end{proof}

\begin{example}
Our count of tree decorations offers a piecewise polynomial version of
the Verlinde formula: if $T$ is any trivalent tree on $n=2k+2$ leaves then the number of $T$-decorations of order $4l$ whose pendent edges have weight~$2l$ is the sum on the right hand side of
Theorem \ref{verlinde}. In particular, if $l=1/2$ then the number of 
$T$-decorations equals $\,2^k = \psi(k,1,\ldots,1)$.
\qed
\end{example}

\section{The zonotopal Cox ring}

In this section we explain the genesis of the present project. Our
point of departure was the work on {\em zonotopal algebra} by Holtz
and Ron \cite{HR} which derives combinatorial formulas for $\psi$
when the linear forms $\ell_j$ and the exponents $u_j$ have the
following special form. Fix a $d \times m$-matrix $C = (c_{ik})$ of
rank $d$ over $K$. Let $H_1,\ldots,H_n$ denote the hyperplanes in
$K^d$ which are spanned by subsets of the columns of $C$. We have $m
\leq n \leq \binom{m}{d-1}$ and the upper bound is attained if the
matrix $C$ is generic. For each hyperplane $H_j$ we fix a nonzero
linear form $\ell_j \in K[z]$ that vanishes on $H_j$. Zonotopal
algebra is concerned with these linear forms $\ell_1,\ldots,\ell_n$
and the ideals $I_u$ generated by certain specific powers of the
$\ell_i$. As before, we write $A= (a_{ij})$ for the $d \times
n$-matrix of coefficients of $\ell_1,\ldots,\ell_n$, and $G \subset
K^n$ is the kernel of $A$.

The Cox-Nagata ring $R^G$ has a special structure which depends on
the given matrix $C$, and if that matrix is generic then $R^G$ will
depend only  on the parameters $d$ and $m$. Also in this special
situation, the ring $R^G$ may fail to be Noetherian. For instance,
this happens when $d=3$ and $m = 9$ because the resulting $n = 36$
linear forms will contain $m=9$ general linear forms.

However, the situation becomes finite and very nice if we restrict
the choice of $u$ to translates of a certain $m$-dimensional
sublattice of $\Z^n$. This sublattice is the image of the following
$m \times n$-matrix ${\bf C}$ with entries in $\{0,1\}$. The
entry of ${\bf C}$ in  row  $k$ and column $j$ is zero if the $k$-th
column of $C$ lies on the hyperplane $H_j$. Let ${\bf e} $ denote
the vector $(1,1,\ldots,1) $ in $\Z^n$. Holtz and Ron \cite{HR} establish
formulas for $\psi(r,{\bf C} v)$, $\psi(r,{\bf C} v - {\bf e})$ and
$\psi(r,{\bf C} v - 2 {\bf e})$. If $C$ is a unimodular matrix then
their formulas are expressed via volumes and lattice points in
zonotopes (whence the term {\em zonotopal algebra}), while in the
general (non-unimodular) case they involve matroid theory. For
instance, $\sum_{r} \psi(r,{\bf C} v )$ is the number of independent
sets in the rank $d$ matroid on $v_1 + \cdots + v_m$ elements
obtained by
 duplicating the $k$-th column of $C$  exactly
$v_k$ times. Likewise, the quantity
 $\sum_{r} \psi(r,{\bf C} v - {\bf e})$  is the number
of bases of that matroid. We shall propose an algebraic explanation
of that result.

We define the {\em zonotopal Cox ring} of the matrix $C$ to be the
subalgebra $Z^G$ of $R^G$ which is the direct sum of all graded
components $R_{r,u}$ where $r \in \N$ and $u$ runs over the lattice
points in the image of ${\bf C}$. We use the notation
$$ Z^G \quad =\,\,\, \bigoplus_{(r,v) \in \Z^{m+1}} \!\!\!\! R^G_{(r,{\bf C} v)} . $$
Clearly, the Cox-Nagata ring $R^G$ is a module over the zonotopal
Cox ring $Z^G$, and for any fixed vector $\omega \in \Z^n$ we can
also consider the submodule
$$ Z^{G,w} \quad =\,\,\, \bigoplus_{(r,v) \in \Z^{m+1}} \!\!\!\! R^G_{(r,{\bf C} v)+w} . $$
We call $Z^{G,w}$ the {\em zonotopal Cox module} of {\em shift} $w$.
Thus the results of Holtz and Ron give formulas for the
$\Z^{m+1}$-graded Hilbert series of the zonotopal Cox ring and the
zonotopal Cox modules with shifts $w = - {\bf e} $ and $w = - 2 {\bf
e}$. Their results have been extended in recent work of Ardila and
Postnikov \cite{AP}.

The generators of the zonotopal Cox ring have the following
description. For each $k \in \{1,2,\ldots,m\}$ let ${\bf c}_k$
denote the $k$-th column of the matrix ${\bf C}$, and let  $\,E_k =
x^{{\bf c}_k}\,$ be the corresponding squarefree monomial. We denote
by   $\, f_k(z) =  \sum_{i=1}^d c_{ik} z_i\,$ the linear form
corresponding to the $k$th column of $C$, and we define $F_k(x,y)$
to be the image of  $f_k(z)$ in the Cox-Nagata ring $R^G$. In other
words, $F_k$ is the numerator of the Laurent polynomial
$\,\sum_{j=1}^n f_k(a_j) \cdot (y_j/x_j) \,$ where $a_j$ is the
$j$th column of the matrix $a_j$. We have
$$ {\rm deg}(E_k) \,\, = \,\, (0,{\bf c}_k) \quad \hbox{and} \quad {\rm deg}(F_k) \,\, = \,\, (1,{\bf c}_k) ,
$$
where ${\bf c}_k$ is the $k$-th column of the matrix ${\bf C}$. In
particular both $E_k$ and $F_k$ lie in the zonotopal Cox ring $Z^G$.
These elements suffice to generate:

\begin{thm}
The zonotopal Cox ring equals $ Z^G  =   K [E_1,\!\ldots \!,  E_m,
F_1,\!\ldots\!,F_m]$.
\end{thm}

\begin{proof}
This lemma is a reinterpretation of the result on exterior
zonotopal algebra in \cite{HR}, which implies that the
$K$-vectorspace $\,(I_{{\bf C}v}^\perp)_r\,$ is spanned by the
products $\,f_1^{h_1} f_2^{h_2} \cdots f_m^{h_m}\,$ where $\,h_1+
h_2 + \cdots + h_m = r\,$ and $h_k \leq v_k$ for all $k$. The image
of that product in $Z^G$ equals $E_1^{v_1-h_1} \cdots E_k^{v_k-h_k}
F_1^{h_1} \cdots F_k^{h_k}$.
\end{proof}

Using constructions from matroid theory, one can select a subset of
the products $\,f_1^{h_1} f_2^{h_2} \cdots f_m^{h_m}\,$ which forms
a basis of the space $\,(I_{{\bf C}v}^\perp)_r \,\simeq\,
Z^G_{(r,{\bf C}{\bf v})}$. The cardinality of that basis is computed
as the value of a {\em multivariate Tutte polynomial}, as explained
in \cite{AP}. In particular, we conclude that the Hilbert function
$(r,v) \mapsto \psi(r,{\bf C}v)$ of the zonotopal Cox ring is a
piecewise polynomial function. We shall now present an explicit formula
for that function.

If $\mu$ is any multiset of positive integers and $s \in \N$ then we
write $\,\Phi(\mu,s) \,$ for the coefficient of $q^s$ in  the
expansion of $\, \prod_{\ell \in \mu} (\sum_{i=0}^{\ell-1} q^i)$.
Thus $\Phi$ is a piecewise polynomial function of degree $|\mu|-1$
in its $|\mu|+1$ arguments. Let $M(C)$ be the rank $d$ matroid on $
\{1,\ldots,m\}$ defined by the matrix $C$. For any $J \subseteq
\{1,\ldots,m\}$ we write ${\rm span}(J) $ for the flat of $M(C)$
spanned by $J$.

\begin{cor} \label{zonoformula}
The Hilbert function of the zonotopal Cox ring $Z^G$ equals
\begin{equation}
\label{zonoformula2} \!  \psi(r,{\bf C}v) \,\,\,\, = \,\,\,\, \sum_I
\Phi\bigl(\,\{ v_i \}_{i \in I}\,,\,\, r \, -  \!\!\! \sum_{j\not\in {\rm span}(I
\cap \{1,2,\ldots,m\})}\!\!\!\!\!\!\!  v_j \, \bigr),
\end{equation}
where the sum is over all independent sets $I$ of the matroid
$M(C)$.
\end{cor}

\begin{proof}
This formula is derived by applying \cite[Theorem 4.2]{HR} to the
rank $d$ matroid on $v_1 +v_2+ \cdots + v_m$ elements which is
obtained from $M(C)$ by duplicating the elements $1,2,\ldots,m$
respectively $v_1,v_2,\ldots,v_m$ times.
\end{proof}

The right hand side of (\ref{zonoformula2}) is obviously a piecewise
polynomial function in $(r,v)$ of degree $d-1$. If we sum this
expression over $r$ from $0$ to $\sum_{i=1}^m v_i$ then the
piecewise nature disappears and we get a polynomial of degree $d$:
\begin{equation}
\label{polydegd}
 \sum_{r \geq 0} \psi(r,{\bf C}v)
\quad = \quad \sum_I \prod_{i \in I} v_i  .
\end{equation}
We now illustrate formula (\ref{zonoformula2}) 
by relating them to earlier examples.

\begin{example}
Let $d=3, m=4$ and $ C \,= \,
\begin{bmatrix}
0 & 0 & 1 & 1 \\
0 & 1 & 0 & 1 \\
1 & 0 & 0 & 1
\end{bmatrix} $.
Then $n=6$ and $A$ is the matrix in Example \ref{sixspecial}. The
zonotopal Cox ring $Z^G$ is the subalgebra of $R^G$ generated by
$L_{124}$, $L_{135}$, $L_{236}$, $L_{456}$, $x_1 x_2 x_3$, $x_1 x_4
x_5$, $x_2 x_4 x_6$ and $x_3 x_5 x_6$. The Hilbert function of $Z^G$ is the following
specialization of $\psi: \N^7 \rightarrow \N$:
\begin{eqnarray*}
& \psi( \,r \, , \,v_1 + v_2, v_1 + v_3, v_2 + v_3, v_1 + v_4, v_2 +
v_4, v_3 + v_4)
\qquad  \\
& =   \,\,\,\, \Phi(\{v_1,v_2,v_3\},r) + \Phi(\{v_1,v_2,v_4\},r {-}
v_3) +
\Phi(\{v_1,v_3,v_4\},r {-} v_2)  \\
 & \quad  + \,\,
\Phi(\{v_2,v_3,v_4\},r  -  v_1)  \,+ \, \sum_{|I| \leq 2} \Phi(
\{v_i \}_{i \in I}, r - \sum_{j \not\in I} v_j ) .
\end{eqnarray*}
If we sum this piecewise quadratic over all $r$ then we get the
cubic polynomial
$$ (\ref{polydegd}) \quad =  \quad
(v_1 +1)(v_2 +1 )(v_3 +1)(v_4 +1) - v_1 v_2 v_3 v_4.
$$
\vskip -.7cm
\qed
\end{example}

A problem that remains open and of interest is to give description
of the Cox-Nagata ring $R^G$ for the configurations considered here.
In algebraic geometry terms, we are concerned with blowing up the
intersection points of a hyperplane arrangement. For instance, we may
ask when $R^G$ is finitely generated, and what is its structure as a module
over the nice subring $Z^G$.

First results in this direction were obtained by Ana-Maria
Castravet. She proved, for instance, that $R^G$ is finitely
generated when $d=4$ and $m=5$, that is, for the blow-up
of $\PP^3$ at the $n=10$ intersection points determined by
five general planes. Her proof is based on the methods
developed in \cite{Cas}.
We close by presenting that example
from our initial perspective in Section~1.

\begin{example} \label{ex:ana}
Consider the system of linear partial differential equations
\begin{eqnarray*}
 I_u \,\,\, = & \bigl\langle
\,\partial_1^{u_1+1}\,, \,\partial_2^{u_2+1}\,,\,\partial_3^{u_3+1}\,, \, \partial_4^{u_4+1}\,,\,
(\partial_1-\partial_2)^{u_5+1}, (\partial_1-\partial_3)^{u_6+1}, \\ & \,\,\,\,\,
 (\partial_1-\partial_4)^{u_7+1},
(\partial_2-\partial_3)^{u_8+1}, (\partial_2-\partial_4)^{u_9+1}, (\partial_3-\partial_4)^{u_{10}+1}\,
\bigr\rangle.
\end{eqnarray*}
We shall compute the number $\psi(r,u) $ of linearly independent polynomial
solutions of degree $r$. To this end, we fix the following matrix
 over  $K = \Q(t)$:
$$ C \quad = \quad
\begin{bmatrix}
\, 0 & 0 & 0 & 1 & t     \, \\
\, 0 & 0 & 1 & 0 & t^2 \, \\
\, 0 & 1 & 0 & 0 & t^3 \, \\
\, 1 & 0 & 0 & 0 & t^4 \,
\end{bmatrix},
$$
and we let $G$ be the corresponding linear
subspace of codimension $4$ in $K^{10}$.
Note that the ten linear forms in $I_u$
arise from the matrix $C$ if we set $t=1$.
A computation with Castravet's generators reveals
 that the subspace $G$ is sagbi. We find that
 the toric algebra ${\rm in}(R^G)$ is generated by
the ten variables
\begin{equation}
\label{xHHH}
\begin{matrix}
123 & 124 & 134 & 234 & 125 & 135 & 235 & 145 & 245 & 345 \\
x_1 & x_2 &  x_3  & x_4  & x_5 & x_6 & x_7 & x_8 & x_9 & x_{10}
\end{matrix}
\end{equation}
which represent the intersection points
$H_j \cap H_k \cap H_l $, and the $15$ monomials
$$
\begin{matrix}
 y_1 x_2 x_3 x_4 &
 y_1 x_5 x_6 x_7 &
y_2 x_5 x_8 x_9 &
 y_3 x_6 x_8 x_{10} &
y_4 x_7 x_9 x_{10}
\\
 y_3 x_4 x_6 x_7 x_8 x_9 &
y_2 x_4 x_5 x_7 x_8 x_{10} &  \! y_1 x_4 x_5 x_6 x_9 x_{10} & \! y_1 x_2 x_3 x_7 x_9 x_{10} &
\! y_2 x_3 x_5 x_6 x_9 x_{10} \\ y_1 x_3 x_5 x_7 x_8 x_{10} &
 y_1 x_2 x_4 x_6 x_8 x_{10} & y_1 x_2 x_6 x_7 x_8 x_9 &
y_1 x_3 x_4 x_5 x_8 x_9 & y_2 x_3 x_4 x_5 x_6 x_7
\end{matrix}
$$
which represent the planes spanned by the ten intersection points in $\PP^3$.
The toric ideal of relations among these $25$ monomials has
$55$ minimal generators, and each of these lifts to a relation in $I^G$.
We compute the $\Z^{11}$-graded Hilbert function of $R^G$ using the technique
explained in Sections 4 and 5, namely by listing
   the facets of the cone $\Gamma$ spanned by the
$25$ monomials. This cone has
$$ f(\Gamma) = (25, \!261, \! 1536, \! 5790, \! 14935, \! 27309, \! 35985,
\! 34247, \! 23276, \! 10989,\! 3419, \! 634, \! 56) $$
and {\tt polymake} supplies  an explicit $14 \times 56$ matrix $M$
such that  $\psi(r,u)$ is the number of integer solutions
$(x,y,z)$ to
$\, (r,u_1,u_2,\ldots,u_{10}, x,y,z) \cdot M \geq 0 $.
\qed
\end{example}

\begin{remark} The $14 \times 56$ matrix $M$ above and other supplementary materials for this
paper are posted at {\tt http://lsec.cc.ac.cn/$\sim$xuzq/cox.html}.
\end{remark}

\bigskip

\noindent {\bf Acknowledgments.}
This research project was carried while both authors visited the
Technical University of Berlin in 2007-08.
 Bernd Sturmfels was supported by an Alexander von Humboldt research prize
and the U.S.~National Science Foundation (DMS-0456960).
Zhiqiang Xu was supported by a Sofia Kovalevskaya prize awarded to Olga Holtz.
We are most grateful to Seth Sullivant and Jenia Tevelev for discussions and comments
which greatly contributed to Section~7.
This article also benefited from interactions with
Klaus Altmann,
Federico Ardila,
Thomas Bauer, Ana-Maria Castravet,
Kristian Ranestad,
Mauricio Velasco and
Volkmar Welker.

\bigskip \medskip 

\noindent {\bf Authors' addresses:}

\medskip

\noindent Bernd Sturmfels, Department of Mathematics,
University of California, \break Berkeley, CA 94720, USA,
{\tt bernd@math.berkeley.edu}

\medskip

\noindent Zhiqiang Xu,
 Academy of Mathematics and System Sciences,
 Chinese Academy of Sciences, Beijing, 100080, China,
 {\tt xuzq@lsec.cc.ac.cn}


\begin{thebibliography}{99}
\bibitem{popov} O.~Popov: The Cox ring of a del Pezzo surface has
rational singularities, {\tt math.AG/0402154}.
\end{thebibliography}
\end{document}